\documentclass[12pt]{amsart}
\usepackage{epsfig}
\usepackage{amsmath}
\usepackage{amsfonts}
\usepackage{amssymb}

\def\C{\mathcal C}
\def\D{\mathcal D}
\def\F{\mathcal F}
\def\K{\mathcal K}
\def\J{\mathcal J}
\def\L{\mathcal L}
\def\H{\mathcal H}
\def\O{\mathcal O}
\def\T{\mathcal T}

\def\1{\mathbf 1}
\def\M{\overline{\mathcal M}}
\def\QQ{\mathbb Q}
\def\ZZ{\mathbb Z}
\def\CC{\mathbb C}
\def\RR{\mathbb R}
\def\Res{\operatorname{Res}}
\def\Ind{\operatorname{Ind}}
\def\hat{\widehat}
\def\tilde{\widetilde}
\def\p{\partial}
\def\a{\alpha}

\def\lan{\langle}
\def\ran{\rangle}
\def\ev{\operatorname{ev}}
\def\ft{\operatorname{ft}}
\def\ct{\operatorname{ct}}
\def\td{\operatorname{td}}
\def\ch{\operatorname{ch}}

\def\Tr{\operatorname{Tr}}

\def\fake{\operatorname{fake}}
\def\tr{\operatorname{tr}}
\def\z{\zeta}

\renewcommand{\Box}{\blacksquare}
\renewcommand{\Delta}{\triangle}

\title[HRR in true Quantum K-theory]
{The Hirzebruch--Riemann--Roch Theorem \\ 
in true genus-0 quantum K-theory}

\author[A. Givental]{Alexander GIVENTAL$^1$}\thanks{$^1$ Partially 
supported by IHES} 
\thanks{$^{1,2}$ This material is based upon work supported by the National 
Science Foundation under Grant DMS-1007164.} 
\author[V. Tonita]{Valentin TONITA$^2$}

\date{June 14, 2011}

\begin{document}

%\maketitle

\begin{abstract}
We completely characterize genus-0 K-theoretic \linebreak 
Gromov--Witten invariants of
a compact complex algebraic manifold in terms of cohomological Gromov--Witten 
invariants of this manifold. This is done by applying (a virtual version of) 
the Kawasaki--Hirzebruch--Riemann--Roch 
formula for expressing holomorphic Euler 
characteristics of orbibundles on moduli spaces of genus-0 stable maps,
analyzing the sophisticated combinatorial structure of inertia stacks of 
such moduli spaces, and employing various quantum Riemann--Roch formulas  
from {\em fake} (i.e. orbifold--ignorant) quantum K-theory of manifold and 
orbifolds (formulas, either previously known from works of Coates--Givental, 
Tseng, and Coates--Corti--Iritani--Tseng, or newly developed for this 
purpose by 
Tonita). The ultimate formulation combines properties of 
overruled Lagrangian cones in symplectic loop spaces (the language, that 
has become traditional in description of generating functions 
of genus-0 Gromov-Witten theory) with a novel framework of {\em adelic 
characterization} of such cones. As an application, we prove that tangent 
spaces of the overruled Lagrangian cones of quantum K-theory carry a natural 
structure of modules over the algebra of finite-difference operators in 
Novikov's variables. As another application, we compute one of such tangent
spaces for each of the complete intersections given by equations of degrees
$l_1,\dots,l_k$ in a complex projective space of dimension 
$\geq l_1^2+\cdots+l_k^2-1$. 

\end{abstract}

\maketitle

\section{Motivation}

K-theoretic Gromov--Witten invariants of a compact complex algebraic manifold
$X$ are defined as holomorphic Euler characteristics of various interesting
vector bundles over moduli spaces of stable maps of compact 
complex curves to $X$. They were first introduced in \cite{GiK}
(albeit, in a limited generality of genus-0 curves mapped to homogeneous
K\"ahler spaces), where it was shown that such invariants 
define on $K^0(X)$ a geometric structure resembling Frobenius structures of
quantum cohomology theory. 

At about the same time, it was shown \cite{GiL} that simplest genus-0 
K-theoretic GW-invariants of the manifold $X$ of complete flags in $\CC^{n+1}$
are governed by the finite-difference analogue of the quantum Toda lattice.
More precisely, a certain generating function for K-theoretic GW-invariants,
dubbed in the literature the {\em J-function} (and depending on $n$ variables,
namely, Novikov's variables $Q_1,\dots,Q_n$ introduced to separate 
contributions of complex curves in $X$ by their degrees)
turns out to be a common eigenfunction (known in representation theory 
as Whittaker's function) of $n$ commuting finite-difference 
operators, originating from the center of the {\em quantized} universal
enveloping algebra $U_qsl_{n+1} $. In quantum cohomology theory, the 
corresponding fact was established by B. Kim \cite{Kim}, who showed that
the cohomological J-function of the flag manifold $X=G/B$ of a complex 
simple Lie
algebra ${\mathfrak g}$ is a Whittaker function of the Langlands-dual Lie 
algebra ${\mathfrak g}'$. The K-theoretic generalization involving quantized 
universal enveloping algebras $U_q{\mathfrak g}'$ was stated in \cite{GiL} 
as a conjecture, and still remains an open problem. 

Foundations for K-theoretic counterpart of GW-theory were laid down by Y.-P. 
Lee \cite{YPLee} in the reasonable generality of arbitrary complex algebraic 
target spaces $X$ (and holomorphic curves of arbitrary genus). While the 
general structure and universal identifies (such as the string equation, or 
topological recursion relations) of quantum cohomology theory carry over 
to case of quantum K-theory, the latter is still lacking certain computational
tools of the former one, and for the following reason. 

The so-called 
{\em divisor equations} in quantum cohomology theory tell that the number of 
holomorphic curves of certain degree $d$ with an additional constraint, that
a certain marked point is to lie on a certain divisor $p$, is 
equivalent to (more precisely, differs by the factor $(p,d)$ from) 
the number of such curves without the marked point and without the constraint. 
Consequently, the dependence of J-functions on Novikov's variables is redundant 
to their behavior as functions on $H^2(X)$ introduced through constraints at
marked points. In particular, differential equations satisfied by the 
J-function in Novikov's variables (e.g. the Toda equations in the case of flag
manifolds) are directly related to the quantum cup-product on $H^*(X)$. 

In K-theory, however, any analogue of the divisor equation seems missing,
and respectively the K-theoretic analogue of the quantum cup-product, and
differential equations of the Frobenius-like structure on $K^0(X)$ are 
completely detached from the way the J-functions depend on Novikov's variables.
Because of this lack of structure with respect to Novikov's variables, it appears
even more perplexing that in examples (such as projective spaces, or flag 
manifolds) the J-functions of quantum K-theory turn out to satisfy interesting 
finite-difference equations.  

\medskip

The idea of computing K-theoretic GW-invariants in cohomological terms is 
naturally motivated by the classical Hirzebruch--Riemann--Roch formula
\cite{Hir}   
\[ \chi (M, V):=\sum_k (-1)^k \dim H^k(M, V) = \int_M \ch(V) \td(T_M).\]
The problem (which is at least a decade old) of putting this idea to work
encounters the following general difficulty. The HRR formula needs to be
applied to the base $M$ which, being a moduli space of stable maps, behaves
as a virtual {\em orbifold} (rather than virtual manifold). The HRR formula
for orbibundles $V$ on orbifolds $M$ was established by Kawasaki \cite{Kaw} and 
expresses the holomorphic Euler characteristic (which is an integer) 
as an integral over the {\em inertia orbifold} $IM$ (rather than $M$ 
itself). The latter is a disjoint union of {\em strata} corresponding to 
points with various types of local symmetry (and $M$ being one of the strata
corresponding to the trivial symmetry). 

When $M$ is a moduli space of stable maps, the strata of the inertia stack 
$IM$ parametrize {\em stable maps with prescribed automorphisms}. It is
reasonable to expect that individual contributions of such strata can be 
expressed
as integrals over moduli spaces of stable maps of quotient curves, and
thus in terms of traditional GW-invariants. However, the mere combinatorics of 
possible symmetries of stable maps appears so complicated (not mentioning
the complexity of the integrands required by Kawasaki's theorem), that 
obtaining a ``quantum HRR formula'' expressing K-theoretic GW-invariants
via cohomological ones didn't seem feasible. 
  
\medskip

In the present paper, we give a complete solution in genus-0 to the problem of
expressing K-theoretic GW-invariants of a compact complex algebraic
manifold in terms of its cohomological GW-invariants. The solution turned
out to be technology-consuming, and we would like to list here those 
developments of the last decade that made it possible.      

One of them is the Quantum HRR formula \cite{CGL, Co} in {\em fake} quantum 
K-theory. One can take the RHS of the classical Hirzebruch--Riemann--Roch 
formula for the definition of $\chi^{\fake}(M,V)$ on an {\em orbifold}
$M$. Applying this idea systematically to moduli spaces of stable maps, one
obtains fake K-theoretic GW-invariants, whose properties are similar to
those of {\rm true} ones, but the values (which are rational, rather
than integer) are different. The formula expresses fake K-theoretic 
GW-invariants in terms of cohomological ones.     
 
Another advance is the Chen--Ruan theory \cite{AGV,ChR} of GW-invariants of 
{\em orbifold target spaces}, and the computation by Jarvis--Kimura 
of such invariants in the case when the target is the quotient of a point
(or more generally a manifold) by the trivial action of a finite group. 
  
Next is the theorem of Tseng \cite{Ts} expressing {\em twisted} 
GW-invariants of orbifold target spaces in terms of untwisted ones. 

Yet, two more ``quantum Riemann--Roch formulas'' of \cite{Co} 
had to be generalized to the case of orbifold targets. This is done 
in \cite{To, To2}.  
 
Finally, our formulation of the Quantum HRR Theorem in true quantum K-theory
is based on a somewhat novel form of describing generating functions of 
GW-theory, which we call {\em adelic characterization}. For a general and 
precise formulation of the theorem, the reader will have to wait until 
Section 6, but here we would like to illustrate the result with an example
that was instrumental in shaping our understanding. 

\medskip 
   
Let 
\[ J=(1-q)\sum_{d\geq 0}\frac{Q^d}{(1-Pq)^n(1-Pq^2)^n\cdots (1-Pq^d)^n}.\]
Here $P$ is unipotent, and stands for the Hopf bundle on $\CC P^{n-1}$, 
satisfying the relation $(1-P)^n=0$ in $K^0(\CC P^{n-1})$.  
It is a power series in Novikov's variable $Q$ with vector coefficients 
which are rational functions of $q$, and take values in $K^0(\CC P^{n-1})$.
It was shown\footnote{Using birational invariance of holomorphic Euler 
characteristcs replacing certain moduli spaces of stable maps 
to $\CC P^{n-1}$ with toric compactifications.} in \cite{GiL} that $J$ 
represents (one value of) the {\em true} 
K-theoretic J-function of $\CC P^{n-1}$.

On the other hand, one can use quantum Riemann--Roch and Lefschetz' theorems
of \cite{Co} and \cite{CCIT} to compute, starting from the cohomological
J-function of $\CC P^{n-1}$, a value of the J-function of the {\em fake} 
quantum K-theory. The result (see Section 10) turns out to be the same: $J$. 
This sounds paradoxical, since --- one can check this directly for $\CC P^1$ 
in low degrees! --- contributions of non-trivial Kawasaki strata neither 
vanish no cancel out. 

In fact this is not a contradiction, for as it turns out, coefficients of
the series $J$ do encode fake K-theoretic GW-invariants, {\em when $J$ is 
expanded into a Laurent series} near the pole $q=1$. Furthermore, when 
$J$ is expanded into a Laurent series near the pole $q=\z^{-1}$, where $\z$ is
a primitive $m$-th root of unity, the coefficients represent certain fake 
K-theoretic GW-invariants of the {\em orbifold target space} $\CC P^{n-1}/\ZZ_m$.
Moreover, according to our main result, these properties altogether completely
characterize those $Q$-series (whose coefficients are vector-valued 
rational functions of $q$) which represent true genus-0 K-theoretic 
GW-invariants of a given target manifold. 

This fact is indeed the result of application
of Kawasaki's HRR formula to moduli spaces of stable maps. Namely, the 
complicated combinatorics of strata of the inertia stacks can be interpreted
as a certain identity which, recursively in degrees, governs the decomposition
of the J-function into the sum of elementary fractions of $q$ with poles at all
roots of unity. The theorem is stated in Section 6 (after the general notations,
properties of quantum K-theory, Kawasaki's HRR formula, and results of fake
quantum K-theory are described in Sections 1--5), and proved in Sections 7 
and 8. 

In Section 10, we develop a technology that allows one to extract concrete 
results from this abstract characterization of quantum K-theory. In particular,
we prove (independently of \cite{GiL}) that the function $J$ is indeed 
the J-function of $\CC P^{n-1}$, as well as similar results for codimension-$k$
complete intersections of degrees $l_1,\dots,l_k$ satisfying 
$l_1^2+\cdots+l_k^2\leq n$. 

Let $q^{Q\p_Q}$ denote the operator of {\em translation} through $\log q$
of the variable $\log Q$. It turns out that for every $s\in \ZZ$, 
\[ \left(Pq^{Q\p_Q}\right)^s J = 
(1-q) P^s\sum_{d\geq 0}\frac{Q^dq^{sd}}{(1-Pq)^n(1-Pq^2)^n\cdots(1-Pq^d)^n}\]
also represent genus-0 K-theoretic GW-invariants of $\CC P^{n-1}$. This example
illustrates a general theorem of Section 9, according to which J-functions of
quantum K-theory are organized into modules over the algebra ${\mathcal D}_q$ 
of finite-difference operators in Novikov's variables. This turns out to be
a consequence of our adelic characterization of quantum K-theory in terms of
quantum cohomology theory, and of the ${\mathcal D}$-module structure (and
hence of the divisor equation) present in quantum cohomology theory.

\section{K-theoretic Gromov--Witten invariants}

Let $X$ be a {\em target} space, which we assume to be a nonsingular
complex projective variety. Let $\M_{g,n}^{X,d}$ denote Kontsevich's moduli 
space of degree-$d$
stable maps to $X$ of complex genus-$g$ curves with $n$ marked points. 
Denote by $L_1,\dots, L_n$ the line (orbi)bundles over $\M_{g,n}^{X,d}$ formed
by the cotangent lines to the curves at the respective marked points.  
When $a_1,\dots, a_n \in K^0(X)$, and $d_1,\dots,d_n\in \ZZ$, we use the
correlator notation
\[ \lan a_1L^{d_1}, \dots, a_n L^{d_n}\ran_{g,n}^{X,d} \]
for the holomorphic Euler characteristic over $\M_{g,n}^{X,d}$ 
of the following sheaf: 
\[ \ev_1^*(a_1) L_1^{d_1} \dots \ev_n^*(a_n)L_n^{d_n} \otimes \O^{vir}.\]
Here $\ev_i:\M_{g,n}^{X,d}\to X$ are the {\em evaluation} maps, and 
$\O^{vir}$ is the {\em virtual structure sheaf} of the moduli spaces of 
stable maps. The sheaf $\O^{vir}$ was introduced by Yuan-Pin Lee \cite{YPLee}. 
It is an element of the Gr\"othendieck group of coherent sheaves on the 
{\em stack} 
$\M_{g,n}^{X,d}$, and plays a role in K-theoretic version of GW-theory of $X$
pretty much similar to the role of the virtual fundamental cycle 
$[\M_{g,n}^{X,d}]^{vir}$ in cohomological GW-theory of $X$. 
According to \cite{YPLee},
the collection of virtual structure sheaves on the spaces $\M_{g,n}^{X,d}$ 
satisfies K-theoretic 
counterparts of Kontsevich--Manin's axioms \cite{KM} for Gromov--Witten 
invariants.  
 
Note that, in contrast with cohomological GW-theory, where the invariants 
are rational numbers, {\em K-theoretic GW-invariants are integers}. 

The following generating function for K-theoretic GW-invariants
is called the {\em genus-0 descendant potential} of $X$:
\[ \F := \sum_{n,d}\frac{Q^d}{n!}\lan t(L),\dots , t(L)\ran_{0,n}^{X,d}.\]
Here $Q^d$ denotes the monomial in the {\em Novikov ring}, the formal  series
completion of the semigroup ring of the {\em Mori cone} of $X$, where 
the monomial represents the degree $d$ of rational curves in $X$, 
and $t$ stands for any Laurent polynomial of one variable, $L$, with vector 
coefficients in $K^0(X)$. 
Thus, $\F$ is a formal function of $t$ with Taylor coefficients in the
Novikov ring.   

\section{The symplectic loop space formalism}   
    
Let $\CC [[Q]]$ be the Novikov ring. Introduce the {\em loop space}
\[ \K := \left[ K^0(X)\otimes \CC (q,q^{-1})\right] \otimes \CC [[Q]].\]
By definition, elements of $K$ are $Q$-series whose coefficients are 
vector-valued rational functions on the complex circle with the coordinate $q$. 
It is a $\CC [[Q]]$-module, but we often suppress Novikov's variables 
in our notation and refer to $\K$ as a linear ``space.'' Moreover, abusing 
notation, we write $\K=K (q,q^{-1})$, where $K=K^0(X)\otimes \CC [[Q]]$. 
We call elements of $\K$ ``rational functions of $q$ with coefficients 
in $K$,'' meaning that they are rational functions in the $Q$-adic sense, i.e. 
modulo any power of the maximal ideal in the Novikov ring.    

We endow $\K$ with {\em symplectic form} $\Omega$, which is a 
$\CC [[Q]]$-valued non-degenerate anti-symmetric bilinear form:
\[ \K\ni f,g \mapsto \ 
\Omega (f,g) = \left[ \Res_{q=0}+\Res_{q=\infty}\right] \ \left( f(q), 
g(q^{-1})\right) \ \frac{dq}{q}.\]
Here $( \cdot , \cdot )$ stands for the K-theoretic intersection pairing on 
$K$:
\[ (a,b):=\chi (X; a\otimes b) = \int_X \td (T_X) \ch (a) \ch (b) .\]
It is immediate to check that the following subspaces in $\K$ are Lagrangian 
and form a {\em Lagrangian polarization}, $\K=\K_{+}\oplus \K_{-}$:
\[ \K_{+}=K[q,q^{-1}],\ \ 
\K_{-}=\left\{ f\in \K\ | \ f(0)\neq \infty,\ f(\infty)=0 \right\} ,\]
i.e. $\K_{+}$ is the space of Laurent polynomials in $q$, and $\K_{-}$ consists
of rational functions vanishing at $q=\infty$ and regular at $q=0$. 

The following generating function for K-theoretic GW-invariants is defined
as a map $\K_{+}\to \K$ and is nick-named the {\em big J-function} of $X$:
\[ \J (t):=(1-q)+t(q)+\sum_a \Phi^a 
\sum_{n,d}\frac{Q^d}{n!}\lan \frac{\Phi_a}{1-qL}, t(L),\dots,t(L)
\ran_{0,n+1}^{X,d}.\]
The first summand, $1-q$, is called the {\em dilaton shift}, the second, 
$t(q)$, the {\em input}, and the sum of the two lies in $\K_{+}$. The remaining
part consists of GW-invariants, with 
$\{ \Phi_a \}$ and $\{ \Phi^a \}$ being any Poincare-dual bases of $K^0(X)$.
It is a formal vector-valued function of $t\in \K_{+}$ with Taylor 
coefficients in $\K_{-}$.
 
Indeed, the moduli space 
$\M_{0,n+1}^{X,d}$ is a ``virtual orbifold'' of finite dimension. In 
particular, in the K-ring of it, the line bundle $L_1$ satisfies 
a polynomial equation, $P(L_1)=0$, with $P(0)\neq 0$.\footnote{On a manifold
of complex dimension $<D$ we would have: $(L-1)^D=0$, i.e. $L$ would be 
unipotent. This may be false on an orbifold, as the minimal polynomial of a 
line bundle can vanish at roots of $1$, but it does not vanish at $0$ since
$L^{-1}$ exists.}
This implies that each correlator tends to $0$ as $q\to \infty$, and therefore
the correlator is a reduced rational function, with the denominator $P(q)$,
and it obviously has no pole at $q=0$.

\medskip

{\tt Proposition.} {\em The big J-function coincides with the differential of
the genus-0 descendant potential, considered as the section of the cotangent
bundle $T^*\K_{+}$ which is identified with the symplectic loop space by 
the Lagrangian polarization $\K=\K_{+}\oplus \K_{-}$ and the dilaton shift 
$f\mapsto f+(1-q)$:}
\[ \J (t) = 1-q + t(q) + d_t \F.\]

{\tt Proof.} To verify the claim, we compute the symplectic inner 
product of the $\K_{-}$-part
of $\J(t)$, with a variation, $\delta t \in \K_{+}$, of the input, and 
show that it is equal to the value of the differential $d_t \F$ on $\delta t$.
Note that, since $\delta t$ has no poles other than $q=0$ or $\infty$, we have 
\begin{align*} \Omega \left(\sum_a \Phi^a \otimes \frac{\Phi_a}{1-qL}, 
\delta t\right) &=   
-\Omega \left(\delta t, \sum_a \Phi^a \otimes \frac{\Phi_a}{1-qL}\right)= \\
- \left[ \Res_{q=0}+\Res_{q=\infty}\right]\ &\ \frac{\sum_a \delta t^a(q)
\Phi_a} {1-q^{-1}L} \frac{dq}{q} 
= \Res_{q=L} \frac{\delta t (q)}{q-L} = \delta t (L).\end{align*}
Therefore the symplectic inner product in question is equal to
\[ \sum_{n,d} \frac{Q^d}{n!} \lan \delta t (L) , t(L), \dots, t(L) 
\ran_{0,n+1}^{X,d} = (d_t\F) (\delta t),\]
as claimed. 
  
\section{Overruled Lagrangian cones}

A Lagrangian variety, $\L$, in the symplectic loop space $(\K, \Omega)$
is called an {\em overruled Lagrangian cone} if $\L$ is a cone with the vertex 
at the origin, and if for every regular point of $\L$, the tangent space, $T$,
is tangent to $\L$ along the whole subspace $(1-q)T$. In particular: (i) 
tangent spaces are invariant with respect to multiplication by $q-1$, (ii) 
the subspaces $(q-1)T$ lie in $\L$ (so that $\L$ is ruled by a 
finite-parametric family of such subspaces), and (iii) the tangent spaces at 
all regular points in a ruling subspace $(q-1)T$ are the same and equal to 
$T$. 

\medskip

{\tt Theorem} (\cite{GiF}).
{\em The range of the big J-function $\J$ of quantum K-theory 
of $X$ is a formal germ at $\J (0)$ of an overruled Lagrangian cone.} 

\medskip

{\tt Proof.} As explained in \cite{GiF}, this is a consequence of the relation
between descendants and {\em ancestors}. 

The {\em ancestor correlators} of quantum K-theory
\[ K^0(X)\ni \tau \ \mapsto \ 
\lan a_1 \bar{L}^{d_1}, \dots, a_n \bar{L}^{d_n}\ran_{g,n}^{X,d} (\tau), \]
are defined as formal power series of holomorphic Euler characteristics
\[ \sum_{l=0}^{\infty} \frac{1}{l!} \chi \left(\M_{g,n+l}^{X,d}; 
\O^{vir}\ev_1^*(a_1)\bar{L}_1^{d_1}\cdots \ev_n^*(a_n)\bar{L}_n^{d_n}
\ev_{n+1}^*(\tau) \cdots \ev_{n+l}^*(\tau) \right),\]
where $\bar{L}_i$, the ``ancestor'' bundles, are pull-backs of the universal 
cotangent line bundles $L_i$ on the Deligne-Mumford space $\M_{g,n}$ by the 
{\em contraction} map $\ct: \M_{g,n+l}^{X,d} \to \M_{g,n}$. The latter map
involves forgetting the map of holomorphic curves to the target space 
as well as the last $l$ marked points. 

The genus-0 {\em ancestor potential} is defined by
\[ \overline{\F}_{\tau} := \sum_{n,d} \frac{Q^d}{n!} \lan \bar{t}(\bar{L}),
\dots, \bar{t}(\bar{L})\ran_{0,n}^{X,d}(\tau)\]
and depends on $\bar{t}$ and $\tau$. The graph of its differential 
is identified in terms of the ancestor version of the big J-function:
\[ \overline{\J}=1-q + \bar{t}(q) + \sum_{a,b} \Phi_a G^{ab}(\tau) \sum_{n,d} 
\frac{Q^d}{n!} \lan \frac{\Phi_b}{1-q\bar{L}}, \bar{t}(\bar{L}),
\dots, \bar{t}(\bar{L})\ran_{0,n+1}^{X,d}(\tau).\]
Here $\left(G^{ab}\right) =\left(G_{ab}\right)^{-1}$, and  
\[ G_{ab}(\tau):=(\Phi_a,\Phi_b)+\sum_{n,d}\frac{Q^d}{n!}
\lan \Phi_a, \tau,\dots, \tau, \Phi_b\ran_{0,2+n}^{X,d}.\]
In the ancestor version of the symplectic loop space formalism, 
the loop space and its polarization $\K=\K_{+}\oplus \K_{-}$ are the same 
as in the theory of descendants, but the symplectic form $\Omega_{\tau}$ 
is based on the pairing tensor $(G_{ab})$ rather than the constant 
Poincare pairing $(\Phi_a,\Phi_b)$. 

Let $\L \subset (\K,\Omega)$ and $\overline{\L}_{\tau}\subset 
(\K,\Omega_{\tau})$ be Lagrangian submanifolds defined by the descendant 
and ancestor J-functions $\J$ and $\overline{\J}$. Then  
\[ \overline{\L}_{\tau} = S_{\tau} \L,\]
where $S_{\tau}: \K \to \K_{\tau}$ is an isomorphism of the symplectic loop 
spaces, defined by the following matrix $S_{\tau}=(S^a_b)$:
\[ S^a_b= \delta^a_b + \sum_{l,d} \frac{Q^d}{l!}
\sum_{\mu} g^{a\mu}\lan \Phi_{\mu}, \tau, \dots, \tau , 
\frac{\Phi_b}{1-qL}\ran_{0,2+n}^{X,d}. \]

It is important that the genus-0 Deligne-Mumford spaces $\M_{0,n}$ are 
manifolds (of dimension $n-3$). Consequently, the line bundles $\bar{L_i}$ are
unipotent. Moreover, at the points $\bar{t}\in \K_{+}$ with $\bar{t}(1)=0$  
the ancestor potential $\overline{\F}_{\tau}$ has all partial 
derivatives of order $<3$ equal to $0$. In geometric terms, the cone 
$\overline{\L}_{\tau}$ is tangent to $\K_{+}$ along $(1-q) \K_{+}$. 
This means that the cone $\L$ is swept by ruling subspaces 
$(1-q)S_{\tau}^{-1}\K_{+}$ parametrized by $\tau \in K$, an that 
each Lagrangian subspace $S_{\tau}\K_{+}$ is tangent to $\L$ along the 
corresponding ruling subspace. The theorem follows. 

\medskip

The proof of the relationship $\overline{\L}=S_{\tau}\L$ is based on 
comparison of the bundles $L_i$ and $\bar{L}_i$, and is quite similar to 
the proof of the corresponding 
cohomological theorem given in Appendix 2 of \cite{CGi}. It uses the 
K-theoretic version of the WDVV-identity introduced in \cite{GiK}, as well
as the {\em string} and {\em dilaton} equations. 

The genus-0 dilaton equation can be derived from the geometric fact
$(\ft_1)_* (1-L_1)=2-n$ about the K-theoretic push-forward along the map 
$\ft_1:\M_{0,n+1}^{X,d}\to \M_{0,n}^{X,d}$ forgetting the first marked point.
It leads to the relation   
\[ \lan 1-L, t(L),\dots,t(L)\ran_{0,n+1}^{X,d} = (2-n)\lan t(L),\dots,t(L)
\ran_{0,n}^{X,d}.\]
The latter translates into the degree-2 homogeneity of $\F$ with respect to
the dilaton-shifted origin, and respectively to the conical property of $\L$.

The string equation is derived from $(\ft_1)_* 1 = 1$ (thanks to 
rationality of the fibers of the forgetting map) and relationships between 
$\ft_1^*(L_i)$ and $L_i$ for $i>1$ (see for instance \cite{GiK}). 
It can be stated as the tangency 
to the cone $\L$ of the linear vector field in $\K$ defined by the operator 
of multiplication by $1/(1-q)$. The operator of multiplication by
\[ \frac{1}{1-q}-\frac{1}{2}=\frac{1}{2}\ \frac{1+q}{1-q}\]
is anti-symmetric with respect to $\Omega$ and thus defines a linear 
{\em Hamiltonian} vector field. Since $\L$ is a cone, this vector field is
also tangent to $\L$, which lies therefore on the zero level of its
quadratic Hamilton function. This gives another, {\em Hamilton-Jacobi form} 
of the string equation.
 
\section{Hirzebruch--Riemann--Roch formula for orbifolds} 

Given a compact complex manifold $M$ equipped with a holomorphic vector 
bundle $E$, the {\em Hirzebruch--Riemann--Roch formula} \cite{Hir} 
provides a cohomological expression for the {\em super-dimension} 
(i.e. Euler characteristic) of the sheaf cohomology:
\[ \chi (M,E):=\dim H^{\bullet}(M,E) = \int_M \td (T_M)\ \ch (E).\]
 The generalization of this formula to the case when $M$ is an orbifold
and $E$ an orbibundle is due to T. Kawasaki \cite{Kaw}. It expresses 
$\chi (M,E)$ as an integral over the {\em inertia orbifold} $IM$ of $M$:
\[ \chi (M,E)=\int_{[IM]} \td (T_{IM}) \ \ch \left( 
\frac{\Tr(E)}{\Tr (\bigwedge^{\bullet} N^*_{IM})} \right).\]

By definition, the structure of an $n$-dimensional complex orbifold on $M$ 
is given by an atlas of local charts $U \to U/G(x)$, the quotients of 
neighborhoods of the origin in $\CC^n$ by (linear) actions of finite {\em local
symmetry groups} (one group $G(x)$ for each point $x\in M$). 

By definition, charts on the inertia orbifold $IM$ have the form
$U^g \to U^g/Z_g(x)$, where $U^g$ is the fixed point locus of $g\in G(x)$, and 
$Z_g(x)$ is the centralizer of $g$ in $G(x)$. For elements $g$ from the same 
conjugacy class, the charts are canonically identified by the action of $G(x)$.
Thus, locally near $x\in M$, connected components of the inertia orbifold 
are labeled by conjugacy classes, $[g]$, in $G(x)$. 
Integration over the fundamental 
class $[IM]$ involves the division by the order of the stabilizer of a typical
point in $U^g$ (and hence by the order of $g$ at least).

Near a point $(x, [g]) \in IM$, the {\em tangent} and {\em normal} orbibundles
$T_{IM}$ and $N_{IM}$ are identified with the tangent bundle to $U^g$ and 
normal bundle to $U^g$ in $U$ respectively.     

The Kawasaki's formula makes use of the obvious lift to $IM$ of the 
orbibundle $E$ on $M$.  By $\bigwedge^{\bullet} N^*_{IM}$, we denoted the 
{\em K-theoretic Euler class} of $N_{IM}$, i.e. the exterior algebra of the 
dual bundle, considered as a $\ZZ_2$-graded bundle (the ``Koszul complex'').

The fiber $F$ of an orbibundle on $IM$ at a point $(x, [g])$ carries the direct
decomposition into the sum of eigenspaces $F_{\lambda}$ of $g$.
By $\Tr (F)$ we denote the {\em trace bundle}\footnote{In fact, 
{\em super-trace}, if the bundle is $\ZZ_2$-graded.},
the virtual orbibundle 
\[  \Tr (F) := \sum_{\lambda} \lambda F_{\lambda} .\]
The denominator in Kawasaki's formula is invertible because $g$ does not have
eigenvalue $1$ on the normal bundle to its fixed point locus. 

Finally, $\td$ and $\ch$ denote the {\em Todd class} and {\em Chern character}.

When $M$ is a global quotient, $\tilde{M}/G$, of a manifold by a finite group,
and $E$ is a $G$-equivariant bundle over $\tilde{M}$, Kawasaki's result 
reduces to Lefschetz' holomorphic fixed point formula for 
super-traces in the sum 
\[ \chi (M,E) = \dim H^{\bullet}(\tilde{M}, E)^G = 
\frac{1}{|G|} \sum_{g\in G} \tr \left(g\ |\ H^{\bullet}(\tilde{M}, E)
\right).\]  

\medskip

The orbifold $M$ is contained in its inertia orbifold $IM$
as the component corresponding to the identity elements of local symmetry 
groups. The corresponding term of Kawasaki's formula is
\[ \chi^{\fake} (M, E):= \int_M \td (T_M) \ \ch (E).\]
We call it the 
{\em fake holomorphic Euler characteristic} of $E$. It is generally
speaking a rational number, while the ``true'' holomorphic Euler characteristic
$\chi (M, E)$ is an integer. 

Note that {\em the right hand side of Kawasaki's formula is the fake 
holomorphic Euler characteristic of an orbibundle, 
$\Tr(E)/\Tr (\bigwedge^{\bullet}(N^*_{IM}))$, on the inertia orbifold.} 

\medskip

Our goal in this paper is to use Kawasaki's formula for expressing 
genus-0 K-theoretic GW-invariants in terms of cohomological ones. 
We refer to \cite{To1} (see also the thesis \cite{To}) for the {\em virtual}
version of Kawasaki's theorem, which justifies application of the formula 
to moduli spaces of stable maps. 

The moduli spaces of stable maps are stacks, i.e. locally are quotients of  
spaces by finite groups. The local symmetry groups $G(x)$ are automorphism 
groups of stable maps. A point in the inertia stack $I\M_{0,n}^{X,d}$ 
is specified by a pair: a stable map to the target space and 
an automorphism of the map. In a sense, a component of the inertia stack 
parametrizes stable maps with prescribed symmetry. 

The components themselves are moduli spaces naturally equipped with 
virtual fundamental cycles and virtual structure sheaves. In fact, 
they are glued from moduli spaces of stable maps of smaller degrees --- 
quotients of symmetric stable maps by the symmetries. Thus the individual 
integrals 
of Kawasaki's formula can be set up as certain invariants of {\em fake} 
quantum K-theory, i.e. fake holomorphic Euler characteristics of certain 
orbibundles on spaces glued from usual moduli spaces of stable maps. 

Our plan is to identify these invariants in terms of conventional ones 
and express them --- and thereby the ``true'' genus-0 K-theoretic 
Gromov-Witten theory --- in terms of cohomological GW-invariants. 

For this, a summary of relevant results about fake quantum K-theory, including 
the Quantum Hirzebruch--Riemann--Roch Theorem of Coates--Givental 
\cite{Co, CGL}, will be necessary.

\section{The fake quantum K-theory}

Fake K-theoretic GW-invariants are defined by
\begin{align*} \lan a_1L^{d_1}, \dots, a_n L^{d_n}\ran_{g,n}^{X,d} &:= \\
\int_{\left[\M_{g,n}^{X,d}\right]^{vir}} & \td \left(T_{\M_{g,n}^{X,d}}
\right)\ \ch \left(  \ev_1^*(a_1) L_1^{d_1} 
\dots \ev_n^*(a_n)L_n^{d_n} \right) ,\end{align*}
i.e. as cohomological GW-invariants involving the Todd class of the 
{\em virtual tangent bundle} to the moduli spaces of stable maps. 

The Chern characters $\ch (L_i)$ are unipotent, and as a result, generating 
function for the fake invariants are defined on the space of formal power 
series of $L-1$. In particular, the big J-function 
\[ \J^{\fake}:=1-q+t(q)+\sum_a \Phi^a\sum_{n,d}\frac{Q^d}{n!}\lan 
\frac{\Phi_a}{1-qL}, t(L),\dots, t(L)\ran_{0,n+1}^{X,d}\]
takes an input $t$\footnote{Note that we still treat our generating functions 
as formal in $t$. In particular, an input here a series in
$q-1$ whose coefficients can be arbitrary as long as they remain ''small''.
In practice they will be the sums of indeterminates (like $t$, which are small 
in their own, $t$-adic topology) with constants taken from the maximal ideal 
of Novikov's ring (and thus small in the $Q$-adic sense).} 
from the space $\K_{+}^{\fake} = K [[q-1]]$ 
of power series in $q-1$ with vector coefficients, and takes values 
in the loop space 
\[ \K^{\fake} := \left\{ \begin{array}{cc} \text{ 
$Q$-series whose coefficients} \\ \text{ are Laurent series in $q-1$} 
\end{array} \right\} .\]
The symplectic form is defined by
\[ \Omega^{\fake}(f,g):=-\Res_{q=1}\left( f(q),g(q^{-1})\right) \frac{dq}{q}.\]
Expand $1/(1-qL)$ into a series of powers of $L-1$:   
\[ \frac{1}{1-qL} = \sum_{k\geq 0} (L-1)^k  
\frac{q^k}{(1-q)^{k+1}}.\]
According to \cite{CGL}, we obtain a Darboux basis: 
\[ \Phi_a (q-1)^k,\ \Phi^a q^k/(1-q)^{k+1},\ a=1,\dots, \dim K^0(X), \ 
k=0,1,2,\dots \]
 Taking
$\K_{-}^{\fake}$ to be spanned over $K$ by $q^k/(1-q)^{k+1}$, we obtain
a Lagrangian polarization of $(\K^{\fake}, \Omega^{\fake})$. As before,
the big J-function coincides, up to the dilaton shift $1-q$, with the graph
of the differential of the genus-0 descendant potential: 
$\J^{\fake}(t)=1-q+t(q)+d_t \F^{\fake}$. 

The range of the function $\J^{\fake}$ forms (a formal germ at $J(0)$ of)
an overruled Lagrangian cone, $\L^{\fake}$.
The proof is based on the relationship \cite{GiF} 
between gravitational descendants and
ancestors of fake quantum K-theory, which looks identical to the one in 
``true'' K-theory (although the values of fake and true GW-invariants 
disagree). 
   
In fact the whole setup for fake GW-invariants can be made purely topological, 
extended to include $K^1(X)$, and moreover, generalized to all 
complex-orientable extraordinary cohomology theories (i.e. complex cobordisms).
In this generality, the quantum Hirzebruch--Riemann--Roch theorem of \cite{Co,
CGL} expresses the fake GW-invariants (of all genera) in terms of the 
cohomological gravitational descendants. The special case we need is stated
below, after a summary of the symplectic loop space formalism of quantum 
cohomology theory. 

\medskip

Take $H=H^{even}(X)\otimes \QQ [[Q]]$, and $(a,b)^H=\int_X ab$. Let
$\H$ denote the space of power $Q$-series whose coefficients are Laurent
series in one indeterminate, $z$. Abusing notation we write: $\H=H((z))$,
(remembering that elements of $\H$ are Laurent series only modulo any power 
of $Q$). Define in $\H$ the symplectic form
\[ \Omega^H(f,g)=\Res_{z=0} \left( f(-z),g(z)\right)^H dz,\]
and Lagrangian polarization
\[ \H_{+}=H[[z]],\ \ \H_{-}=z^{-1} H[z^{-1}] .\]
Using Poincare-dual bases of $H$, and the notation $\psi = c_1(L)$, 
we define the big J-function of cohomological GW-theory 
\[ \J^H = -z+t(z)+\sum_a \phi^a \sum_{n,d} \frac{Q^d}{n!}
\lan \frac{\phi_a}{-z-\psi}, t(\psi),
\dots, t(\psi)\ran_{0,n+1}^{X,d}\]
It takes inputs $t$ from $\H_{+}$, takes values\footnote{The previous footnote 
about fake K-theory applies here too. In particular, for the formal function, 
to assume {\em values} in $\H$ merely means that the {\em coefficients} of it 
as a formal $t$-series become Laurent series in $z$ when reduced modulo 
a power of $Q$.} in $\H$,  
and coincides with the graph of differential of the cohomological genus-0 
descendant potential, $\F^{H}$, subject to the dilaton shift $-z$:
$\J^H(t) =-z-t(z)+d_t \F^H$.  Here
\[ \F^H:=\sum_{n,d}\frac{Q^d}{n!} \lan t(\psi),\dots, t(\psi)\ran_{0,n}^{X,d},
\]
where for $a_i\in H$ and $d_i\geq 0$, we have:
\[ \lan a_1 \psi^{d_1},\dots, a_n \psi_n^{d_n} \ran_{0,n}^{X,d}:=
 \int_{\left[\M_{0,n+1}^{X,d}\right]} \ev_1^*(a_1) c_1(L_1)^{d_1} \cdots 
\ev_n^*(a_n) c_1(L_n)^{d_n} .\]
The range of the function $\J^H$ is a Lagrangian cone, $\L^H\subset \H$,
overruled in the sense that its tangent spaces, $T$, are tangent to $\L^H$
along $zT$ (see Appendix 2 in \cite{CGi}). 

%The following HRR formula in fake quantum K-theory makes use of the ordinary
%cup-product operation in $H$.   

\medskip

{\tt Theorem} (\cite{CGL}, see details in \cite{Co}).
{\em  Denote by $\Delta$ the 
Euler--Maclaurin asymptotic of the infinite product
\[ \Delta \sim \prod_{\text{Chern roots $x$ of $T_X$}} \ 
\prod_{r=1}^{\infty} \frac{x-rz}{1-e^{-x+rz}} .\]
Identify $\K^{fake}$ with $\H$ using the Chern character 
isomorphism $\ch: K\to H$ and $\ch (q)=e^z$.
Then $\L^{\fake}$ is obtained 
from $\L^H$ by the pointwise multiplication on $\H$ by $\Delta$:
\[ \ch \left( \L^{\fake} \right) = \Delta \L^H.\]}

\vspace{-3mm}
{\tt Remarks.} (1) Given a function $x\mapsto s(x)$, the {\em Euler--Maclaurin 
asymptotics} of $\prod_{r=1}^{\infty} e^{s(x-rz)}$ is obtained by the 
formal procedure:
\begin{align*} \sum_{r=1}^{\infty} s(x-rz)&=\left(\sum_{r=1}^{\infty} 
e^{-rz\p_x} \right) s(x) = \frac{z\p_x}{e^{z\p_x}-1} (z\p_x)^{-1} s(x) \\ 
& =\frac{s^{(-1)}(x)}{z}
-\frac{s(x)}{2}+\sum_{k=1}^{\infty} \frac{B_{2k}}{(2k)!} s^{(2k-1)}(x) 
z^{2k-1}, \end{align*}
where $s^{(k)}=d^ks/dx^k$, $s^{(-1)}$ is the anti-derivative 
$\int_0^x s(\xi) d\xi$, and $B_{2k}$ are Bernoulli numbers.
Taking $e^{s(x)}$ to be the Todd series, $x/(1-e^{-x})$,
and summing over the Chern roots $x$ of the tangent bundle $T_X$, we get:
\[ \Delta = \frac{1}{\sqrt{\td (T_X)}} 
\exp \left\{ \sum_{k\geq 0}\sum_{l\geq 0} s_{2k-1+l} 
\frac{B_{2k}}{(2k)!} \ch_l(T_X) z^{2k-1} \right\},\]
where the coefficients $s_l$ hide another occurrence of Bernoulli numbers:
\[ \text{\Large $e^{\sum_{l\geq 0} s_l x^l/l!}$} = 
\frac{x}{1-e^{-x}}= 
1+\frac{x}{2}+\sum_{l=1}^{\infty} \frac{B_{2l}}{(2l)!} x^{2l}.\]

(2) Note that neither $\ch: \K\to \H$ nor $\Delta : \H \to \H$ is symplectic:
the former because $(a,b)^{\fake}=(\ch(a), \td (T_X) \ch(b))^H\neq 
(\ch(a),\ch(b))^H$,
the latter because of the factor $\td (T_X)^{-1/2}$. However the composition 
$\ch^{-1}\circ \Delta: \H\to \K^{\fake}$ is symplectic.

(3) The transformation between cohomological and K-theoretic J-functions 
(or descendant potentials) encrypted by the theorem, involves three 
aspects. One is the transformation $\Delta$, while the other
two are the changes of the polarization and dilaton shift. Namely, 
$\ch^{-1}: \H\to \K^{\fake}$ maps $\H_{+}$ to $\K_{+}$ but does {\em not}
map $\H_{-}$ to $\K_{-}^{\fake}$, and there is 
a discrepancy between the dilaton shifts: $\ch^{-1}(-z)=\log q^{-1} \neq 1-q$.
 
(4) Since $\L^{\fake}$ is an overruled cone, it is invariant under the 
multiplication by the ratio $(1-q)/\log q^{-1}$. This shows one way of 
correcting for the discrepancy in dilaton shifts.

(5) The proof of the theorem does not exploit any properties of overruled 
cones. One uses the family $\td_{\epsilon}(x):=\epsilon x/(1-e^{-\epsilon x})$ 
of ``extraordinary'' Todd classes to interpolate between cohomology and 
K-theory, and establishes an infinitesimal version of the theorem. For this,
the {\em twisting classes} $\td_{\epsilon} (T_{\M_{g,n}^{X,d}})$ of the moduli
spaces are expressed in terms of the descendant classes by applying the 
Gr\"othendieck--Riemann--Roch formula to the fibrations 
$\ft_{n+1}: \M_{g,n+1}^{X,d}\to \M_{g,n}^{X,d}$. 

We refer for all details to the dissertation \cite{Co}. However, in Section 8, 
we indicate geometric origins of the three changes described by the 
theorem: the change in the position of the cone, 
in the dilaton shift, and in the polarization.   

\section{Adelic characterization of quantum K-theory}

Recall that point $f\in \K$ is a series in the Novikov variables, $Q$, 
with vector coefficients which are rational functions of $q^{\pm 1}$. 
For each $\z\neq 0,\infty$, we expand (coefficients of) $f$ in a Laurent
series in $1-q\zeta$ and thus obtain the {\em localization} $f_{\z}$ near
$q=\zeta^{-1}$. Note that for $\z=1$, the localization lies in the loop 
space $\K^{\fake}$ of fake quantum K-theory. The main result of the present
paper is the following theorem, which provides a complete characterization 
of the true quantum K-theory in terms of the fake one. 

\medskip

{\tt Theorem.} {\em 
Let $\L \subset \K$ be the overruled Lagrangian cone of quantum
K-theory of a target space $X$. If $f\in \L$, then
the following conditions are satisfied:

{\em (i)} 
$f$ has no pole at $q=\zeta^{-1}\neq 0,\infty$ unless $\zeta$ is a root
of $1$.    

{\em (ii)} When $\z=1$, the localization $f_{\z}$ lies in 
$\L^{\fake}$. 

In particular, the localization $\J(0)_1$ at $\z=1$ of the value of the 
J-function with the input $t=0$ lies in $\L^{\fake}$. In the tangent space to 
$\L^{\fake}$ at the point $\J(0)_1$, make the change $q\mapsto q^m$, 
$Q^d\mapsto Q^{md}$, and denote by $\T$ the resulting subspace in 
$\K^{\fake}$.  
Let $\nabla_{\zeta}$ denote the Euler--Maclaurin asymptotics as $q\z\to 1$
of the infinite product:
\[ \nabla_{\z}\ \sim_{q\z \to 1} \ 
\prod_{\stackrel{\text{\footnotesize \em K-theoretic Chern}}
{\text{\footnotesize \em roots $P$ of $T^*_X$}}}
 \frac{\prod_{r=1}^{\infty} (1-q^{mr}P)}{\prod_{r=1}^{\infty}(1-q^rP)}.\]

{\em (iii)} If $\z\neq 1$ be a primitive $m$-th root of 1, then  
$\left(\nabla_{\z}^{-1} f_{\z}\right) (q/\z) \in \T$.

\medskip

Conversely, if $f\in \K$ satisfies conditions {\em (i),(ii),(iii),} then 
$f\in \L$.}

\medskip

{\tt Remarks.}
(1) The cone $\L$ is a formal germ at $\J(0)$. The statements (direct and
converse) about ``points'' $f\in \L$ are to be interpreted in the spirit of
formal geometry: as statements about {\em families} based at $\J(0)$. 

(2) K-theoretic Chern roots $P$ are characterized by 
$\ch (P)=e^{-x}$ where $x$ are cohomological Chern roots of $T_X$.  

(3) After the substitution $q\zeta =  e^z$ the infinite product becomes
\[ \prod_{\text{Chern roots $x$}} \prod_{k=1}^{m-1}
\prod_{r=0}^{\infty} (1-\z^{-k}e^{kz} e^{-x+mrz})^{-1}.\]
The Euler--Maclaurin expansion has the form
\[ \log \nabla_{\z} = \frac{s^{(-1)}}{mz} + \frac{s}{2} +
\sum_{k>0} \frac{B_{2k}}{(2k)!} (mz)^{2k-1} s^{(2k-1)},\]
where $s$ also depends on $z$ as a parameter:
\[ s (x,z) = - \log \prod_x \prod_{k=1}^{m-1} (1-\z^{-k}e^{kz}e^{-x}).\]
Note that since $x$ are nilpotent, $s(x,z)$ is polynomial in $x$ with
coefficients which expand into power series of $z$. The scalar factor
of $\nabla_{\z}$ is $e^{s(0,0)/2} = m^{-\dim X/2}$ since for each of $\dim X$
Chern roots,
\[ \lim_{x\to 0}\prod_{k=1}^{m-1} (1-\z^{-k}e^{-x}) 
=\lim_{x\to 0}\frac{1-e^{-mx}}{1-e^{-x}} = m.\]

(4) The (admittedly clumsy) definition of subspace $\T$ can be clarified as
follows. The tangent space to $\L^{\fake}$ at the point $\J(0)_1$ is the range
of the linear map $S^{-1}: \K_{+}^{\fake} \to \K^{\fake}$, where
$S^{-1}$ is a matrix Laurent series in $q-1$ with coefficients in the Novikov
ring (see Section 3). Let $\tilde{S}$ be obtained from $S$ by the change
$q\mapsto q^m$, $Q^d\mapsto Q^{md}$.
Then $\T := \tilde{S}^{-1}\K^{\fake}_{+}$. 
   
(5) The condition (iii) seems ineffective, since it refers to a tangent 
space to the cone $\L^{\fake}$ at a yet unknown point $\J(0)_1$. However,
we will see later that the three conditions together allow one, at least in 
principle, to compute the values $\J(t)$ for any input $t$, assuming that the
cone $\L^{\fake}$ is known, in a procedure recursive on degrees of stable maps.
 In particular, this applies to $J(0)_1$. The cone $\L^{\fake}$,
in its turn, is expressed through $\L^H$, thanks to the quantum HRR theorem 
of the previous section, by a procedure which in principle has a similar
recursive nature. Altogether, our theorem expresses all genus-0 
K-theoretic gravitational descendants in terms of the cohomological ones.
Thus this result indeed qualifies for the name:
{\em the Hirzebruch--Riemann--Roch theorem of true genus-0 quantum K-theory.}

\medskip

We describe here a more geometric (and more abstract) formulation of 
the theorem using the {\em adelic} version of the symplectic loop space 
formalism. 

For each $\z\neq 0,\infty$, let $\K^{\z}$ be the space of
power $Q$-series with vector Laurent series in $1-q\z$ as coefficients.
Define the symplectic form
\[ \Omega^{\z}(f,g)=-\Res_{q=\z^{-1}} \left( f(q),g(q^{-1})\right) 
\ \frac{dq}{q},\]
and put $\K_{+}^{\z}:=K [[1-q\z]]$. 
The {\em adele space} is defined as the subset in the Cartesian product:
\[ \hat{\K}\subset \prod_{\z\neq 0,\infty} \K^{\z}\]
consisting of collections $f_{\z}\in \K^{\z}$ such that, modulo any power 
of Novikov's variables, $f_{\z}\in \K_{+}$ for all but finitely many 
values of $\z$. The adele space is equipped with the product symplectic form:
\[ \hat{\Omega}(f,g)=-\sum \Res_{q=\z^{-1}} \left( f_{\z}(q),g_{\z}(q^{-1})
\right)\ \frac{dq}{q}.\]

Next, there is a map $\K \to \hat{\K}: f\mapsto \hat{f}$, which to a rational 
function of $q^{\pm 1}$ assigns the collection $(f_{\z})$ of its 
localizations at $q=\z^{-1}\neq 0,\infty$. Due to the residue theorem, 
the map is symplectic:
\[ \Omega (f,g) = \hat{\Omega} (\hat{f},\hat{g}).\] 
Given a collection $\L^{\z}\subset (\K^{\z},\Omega^{\z})$ of overruled 
Lagrangian cones such that modulo any power of Novikov's variables, 
$\L^{\z}=\K_{+}^{\z}$ for all but finitely many values of $\z$, the product
$\prod_{\z\neq 0,\infty}\L^{\z}\subset \hat{\K}$ 
becomes an {\em adelic} overruled 
Lagrangian cone in the adele symplectic space. 

In fact, ``overruled'' implies invariance of tangent spaces under 
multiplication by $1-q$. Since $1-q$ is invertible at $q=\z^{-1}\neq 1$,
all $\L^{\z}$ with $\z\neq 1$ must be linear subspaces.

According to the theorem, {\em the image $\hat{\L}\subset \hat{\K}$ 
of the cone 
$\L\subset \K$ under the map $\hat{\ }: \K\to \hat{\K}$ followed by a suitable 
adelic (pointwise) completion, is an adelic overruled Lagrangian cone:
\[ \hat{\L} = \prod_{\z\neq 0,\infty} \L^{\z},\]
where $\L^{\z}=\K^{\z}_{+}$ unless $\z$ is a root of $1$, $\L^{\z}=\L^{\fake}$
when $\z=1$, and $\L^{\z}=\nabla_{\z} \T^{\z}$ when $\z\neq 1$ is a root of 
$1$, $\T^{\z}\subset \K^{\z}$ being obtained from the subspace $\T\subset
\K^{\fake}$ (described in the theorem) by the isomorphism 
$\K^{\fake} \to \K^{\z}$ induced by the change $q\mapsto q\z$.}
 
\medskip

{\tt Corollary.} {\em Two points $f,g\in \L$ lie in the same ruling space of 
$\L$ if and only if their expansions $f_1,g_1$ near $q=1$ lie in the same 
ruling space of $\L^{\fake}$.}

\medskip

{\tt Proof.} If $f_1,g_1$ lie in the same ruling space of $\L^{\fake}$, then
$\epsilon \hat{f}+(1-\epsilon)\hat{g} \in \hat{L}$ for each value of $\epsilon$,
and therefore, by the theorem, the whole line $\epsilon f+(1-\epsilon) g$ 
lies in $\L$. The converse is, of course, also true: if the line through 
$f,g$ lies in $\L$ then the line through $f_1,g_1$ lies in $\L^{\fake}$. 
It remains to notice that ruling spaces of $\L$ or $\L^{\fake}$ are 
{\em maximal} linear subspaces of these cones (because this is true modulo 
Novikov's variables, i.e. in the classical K-theory). $\Box$

\section{Applying Kawasaki's formula}

Here we begin our proof of the theorem formulated in the previous section.
The big J-function (see Section 2) consists of the {\em dilaton shift} $1-q$, 
the {\em input} $t(q)$, and holomorphic Euler characteristics of bundles
on virtual orbifolds $\M_{0,n+1}^{X,d}$. The Euler characteristics can be 
expressed,
by applying Kawasaki's formula, as sums of fake holomorphic Euler 
characteristics over various strata of the inertia stacks 
$I\M_{0,n+1}^{X,d}$. A point in the inertia stack is represented by 
a stable map with symmetry (an automorphism, possibly trivial one). 
A stratum is singled out by the combinatorics of such a curve with symmetry. 
Figure 1 below is our book-keeping device for cataloging all the strata.  
  
\begin{figure}[htb]
\begin{center}
\epsfig{file=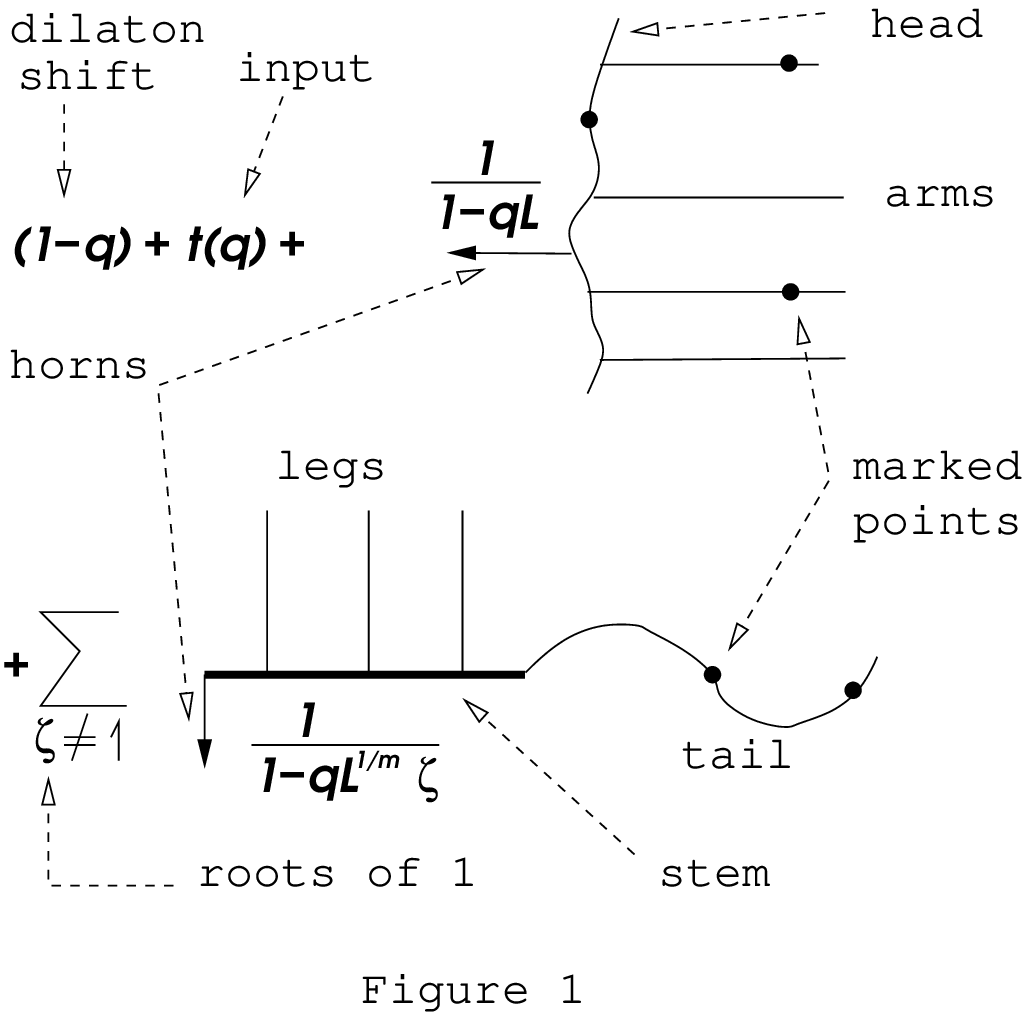}
\end{center}
\end{figure}

Let's call what's written in a given seat of a correlator the {\em content} of 
that seat. In the J-function, the content of the first marked point
has the factor $1/(1-qL)$. We call this marked point the {\em horn}. 

Given a stable map with symmetry, we focus out attention on the horn. The 
symmetry preserves the marked point and acts on the {\em co}tangent line 
at this point with an eigenvalue, which we denote $\zeta$. In Figure 1,
contributions of strata with $\zeta=1$ are separated from those where 
$\zeta\neq 1$, in which case $\zeta$ is a primitive root of 1 of certain order 
$m\neq 1$. 

When $\zeta=1$, the symmetry is trivial on the irreducible component of the
curve carrying the horn. In the curve, we single out the {\em maximal} 
connected subcurve containing the horn on which the symmetry is trivial, 
and call this subcurve (and the restriction to it of the stable map) 
the {\em head}. 

The heads themselves are stable maps without symmetry, and are 
parametrized by moduli spaces $\M_{0,n'+1}^{X,d'}$. Apart from the horn,
the $n'$ marked points are either marked points of the whole curve
or the nodes where ``arms'' are attached. An {\em arm} a stable map 
obtained as a connected component of what's left of the original curve
when the head is removed. The arm has its own horn ---  the nodal point
where it is attached to the head. An arm can be {\em any} stable map with
symmetry, with the only restriction: at its horn, the eigenvalue of the 
symmetry $\neq 1$ (because otherwise the head could be increased).    

In Figure 1, contributions of the strata with the eigenvalues $\zeta\neq 1$
are appended into the sum $\sum$. If $g$ denotes the symmetry of the stable 
map, and $\zeta$ is a primitive $m$th root of $1$, then $g^m$ acts trivially
on the component carrying the horn. We single out the {\em maximal} connected 
subcurve on which $g^m$ is trivial. 
Then the restriction of the stable map to it
has $g$ as a symmetry of order $m$. We call the quotient stable map (of the
quotient curve) the {\em stem}. We'll come back soon to a detailed 
discussion of ``legs'' and ``tails'' attached to the stems. 

Let us denote by $L$ 
the universal cotangent line (on the moduli space of stems) at the horn.
The content in the fake holomorphic Euler characteristic represented by this
term in the sum $\sum$ has the factor $1/(1-qL^{1/m}\zeta)$. Indeed, if $L'$ 
denotes the universal cotangent line to the original stable map, restricted
to the stratum in question, then (in the notation of Kawasaki's formula
in Section 4) $\ch (\Tr L') = \zeta e^{c_1(L)/m}$.

Note that $c_1(L)$ is nilpotent on each of the stem or head spaces. 
Thus, Figure 1 provides the decomposition of $\J$
into the Laurent polynomial part $1-q+t(q)$ and elementary fractions 
$1/(1-q\zeta)^r$ at different poles $q=\zeta^{-1}\neq 0,\infty$.  
We are ready for our first conclusion.

\medskip

{\tt Proposition 1.} {\em The localization $\J_{\z}$ at $\z=1$ lies in the 
cone $\L^{\fake}$ of fake quantum K-theory.} 
 
\medskip

{\tt Proof.} Denote by $\tilde{t}(q)$ the sum of $t(q)$ and of all the terms
of $\sum$ with $\z\neq 1$. Note that in genus $0$, stable maps of degree $0$ 
have no non-trivial automorphisms. So all terms of the sum $\sum$ have 
non-zero degrees. This shows that (thanks to Novikov's variables) the whole 
sum $\tilde{t}$ makes sense as a $q$-series lying in $\K^{\fake}_{+}$, and is
 ``small'' in the $(t,Q)$-adic sense, hence qualifying on the role of an 
input of fake quantum K-theory of $X$. We claim that the whole sum shown on 
Figure 1 is the value of J-function of this fake theory with the input 
$\tilde{t}$.   

Indeed, examine contributions into the virtual Kawasaki formula \cite{To1}
of the terms with $\z=1$. 
Denote by $L_{-}$ the cotangent line at a marked point of the head. 
When the marked point of the head is that of the original curve, the content 
of it is $t(L_{-})$. When this is a node where an arm is attached, 
denote by $L_{+}$ 
the cotangent line to the arm. The only ingredients that do not factor 
into separate contributions of the head and of the arms are
\[ \frac{\sum_a\Phi_a\otimes \Phi^a}{1-L_{-}\Tr(L_{+}')}. \]
The top comes from the gluing of the arm to the head, 
and the bottom from the smoothing of the curve at the node, as a mode of 
perturbation {\em normal} to the stratum of the inertia orbifold.
We conclude that the content of the marked point of the head correlator
is exactly $\tilde{t}(L_{-})$. 

Thus $\J(t)$ is represented as 
\[ 1-q + \tilde{t}(q)+\sum_a \Phi^a \sum_{n',d'}\frac{Q^{d'}}{n'!} 
\lan \frac{\Phi^a}{1-qL}, \tilde{t}(L),\dots, \tilde{t}(L)
\ran_{0,n'+1}^{X, d'} = \J^{\fake}(\tilde{t}),\]
since the correlators come from the fake K-theory of $X$. $\Box$  

\medskip

Let us return to the term with $\z\neq 1$. The stem curve itself is typically
the quotient of $\CC P^1$ by the rotation through $\z$ about two points: the
horn and one more --- let's call it the {\em butt} --- where the eigenvalue
of the symmetry on the cotangent line is $\z^{-1}$. In fact the stem can 
degenerate into the quotient of a {\em chain} of several copies of $\CC P^1$,
with the same action of the symmetry on each of them, and connected 
``butt-to-horn'' to each other (and even further, with other irreducible 
components attached on the ``side'' of the chain, see Figure 2 in the next 
section).  In this case the butt of the stem is that of
the last component of the chain. 
The butt can be a regular point of the whole curve, a 
marked point of it, or a node where the {\em tail} 
is attached (see Figure 1). The tail can be any stable map with any symmetry,
except that at the point where it is attached to the 
stem, the eigenvalue of the symmetry cannot be equal to $\z$. (Otherwise the 
stem chain could be prolonged.) In Figure 1, put 
${\delta t}(q)=1-q+t(q)+\check{t}(q)$, where $\check{t}(q)$ is
the sum of all remaining terms except the one with the pole
at $q=\z^{-1}$ (with this particular value of $\z$). We claim that 
the expansion $\J_{\z}$ of the big J-function near $q=\z^{-1}$ has the form 
\[ \delta{t}(q)+
\sum_a\Phi_a\sum_{n,d}\frac{Q^{md}}{n!}\left[ \frac{\Phi^a}
{1-qL^{1/m}\z}, T(L),\dots,T(L),\delta t(L^{1/m}/\z)\right]_{0,n+2}^{X,d},\]
where $\left[\dots\right]$ are certain correlators of ``stem'' theory, 
and $T(L)$ are leg contributions, both yet to be identified.

Indeed, let $L_{+}$ denote the cotangent line at the butt of the stem, 
and $L_{+}'$ its counterpart on the $m$-fold cover. When the butt is a marked 
point, its content is $t(L_{+}^{1/m}\z)$, and when it is the node with a tail 
attached, then it is $\check{t}(L_{+}^{1/m}\z)$. This is because
$\ch (\Tr L_{+}') =\z e^{c_1(L_{+})/m}$. The case when the butt is a regular 
point on the original curve but a marked point on the stem, can be compared to
the case when the butt is a marked point on the original curve as well. 
In the former case, the conormal bundle to the stratum of stable maps 
with symmetry is missing, comparing to the latter case, the line $L_{+}'$.  
In other words, one can replace the former contribution with the latter one,
by taking the content at the butt to be $1-L_{+}^{1/m}\z$, i.e. the 
K-theoretic Euler factor corresponding to the conormal line bundle $L_{+}'$. 
We summarize our findings.

\medskip

{\tt Proposition 2.} {\em The expansion $\J_{\z}$ of $\J$ near $q=\z^{-1}$
is a tangent vector to the range of the fake J-function of the ``stem'' 
theory at the ``leg'' point, $T$.}

\medskip
    
Our next goal is to understand leg contributions $T(L)$. 

\medskip

{\tt Proposition 3.} {\em Let $\tilde{T}(L)$ denote the arm contribution 
$\tilde{t}(L)$ computed when the input $t=0$. Then
\[ T(L)= \Psi^m \left(\tilde{T}(L)\right).\]}

%\medskip

We remind that Adams' operations $\Psi^m$ are additive and multiplicative
endomorphisms of K-theory acting on a line bundle by $\Psi^m(L)=L^m$.
In this proposition, $\Psi^m$ acts not only on $L$ and elements of $K^0(X)$, 
but also by $\Psi^m(Q^d)=Q^{md}$ on Novikov's variables. 

\medskip

{\tt Proof.} The {\em legs} of a stable map with an automorphism, $g$, 
of order $m\neq 1$ on the cotangent line at the horn, are obtained by 
removing the stem (and the tail). Each leg shown in Figure 1 represents
$m$ copies of the same stable map glued to the $m$-fold cover of the stem
and cyclically permuted by $g$. The automorphism $g^m$ preserves each copy 
of the leg but acts non-trivially on the cotangent line at the horn of 
the leg (i.e. the point of gluing), since otherwise the stem could be 
extended. The only other restriction on what a leg could be is that it cannot
carry (or be) a marked point of the original curve, since the numbering of the 
$m$ copies of the marked point would break the symmetry. This identifies each 
of the $m$ copies of a leg with an arm carrying no marked points. 
       
As in the proof of Proposition 1, denote by $L_{-}$ and $L'_{+}$ the 
cotangent 
lines at the point of gluing to the $m$-fold cover of the stem and to the leg 
respectively. Then the smoothing perturbation at the node of the curve 
with symmetry represents a direction normal to the stratum of symmetric curves.
In the denominator of the virtual Kawasaki formula \cite{To1}, 
it is represented by one 
Euler factor $1-L_{-}\Tr (L'_{+})$ for each copy of the leg. As in the case 
of arms, the gluing factor has the form
\[ \frac{\sum_a \Phi_a\otimes \Phi^a}{1-L_{-}\Tr (L'_{+})}.\]
Then $\ch(\Phi^a)$ and $\ch\left(\Tr(L'_{+})\right)$ are integrated out 
over the moduli space of legs,
and the leg contributes into the fake Euler characteristics over the space
of stems through $\Phi_a$ and $L_{-}$. We claim however that the contribution
of the gluing factor into the stem correlator has the form
\[ \frac{\Psi^m(\Phi_a)\otimes \Phi^a}{1-L_{-}^m\Tr(L'_{+})}.\]
This follows from the following general lemma.

\medskip

{\tt Lemma.} {\em Let $V$ be a vector bundle, and $g$ the automorphism of
$V^{\otimes m}$ acting by the cyclic permutation of the factors. Then
\[ \Tr (g\ |\ V^{\otimes m}) = \Psi^m(V).\]}

\medskip
         
We conclude that the contribution of the leg into stem correlators
is obtained from $\tilde{t}(L_{-})$ (the contribution of the arm into head 
correlators) by computing it at the input $t=0$ (this eliminates those arms
that carry marked points), then applying $\Psi^m$, and also replacing $Q^d$
with $Q^{md}$, because the total degree of the $m$ copies of a leg is $m$ times
the degree of each copy. $\Box$

\medskip

{\tt Proof of Lemma.} It suffices to prove it for the universal $U_N$-bundle,
or equivalently, for the vector representation $V=\CC^N$ of $U_N$.  
Computing the value at $h\in U_N$ of the character of $\Tr (g\ |\ V^{\otimes 
m})$, considered as a representation of $U_N$, we find that it is equal to 
$\tr (g h^{\otimes m})$, because $g$ and $h^{\otimes m}$ commute.
Let $e_i$ denote eigenvectors of $h$ with eigenvalues $x_i$.
A column of the matrix of $gh^{\otimes m}$ in the basis 
$e_{i_1}\otimes \cdots \otimes e_{i_N}$ has zero diagonal entry unless 
$i_1=\cdots=i_N$. 
Thus, $\tr (g h^{\otimes N})=
x_1^m+\cdots+x_N^m$. This is the same as the trace of $h$ on $\Psi^m(V)$.
$\Box$

\medskip
 
{\tt Remark.} The lemma can be taken for the definition of Adams' operations. 
For a permutation $g$ with $r$ cycles of lengths $m_1,\dots,m_r$, it implies: 
\[ \Tr (g\ |\ V^{\otimes m}) = \Psi^{m_1}(V)\otimes \dots \otimes 
\Psi^{m_r}(V).\]  

\medskip

{\tt Proposition 4.}
{\em Propositions 1,2,3 unambiguously determine the big J-function $\J$ in 
terms of stem and head correlators.} 
 
\medskip

{\tt Proof.} Figure 1 can be viewed as a recursion relation that reconstructs 
$J(t)$ by induction on degrees $d$ of Novikov's monomials $Q^d$ (in the sense 
of the natural partial ordering on the Mori cone). 
The key fact is that in genus $0$, 
constant stable maps have no non-trivial automorphisms (and have $>2$ marked 
or singular points). Consequently, arms which are not marked points, as well as
legs, or stems with no legs attached, must have non-zero degrees. As a result,
setting $t=0$, one can reconstruct $\J(0)$ up to degree $d$ from head and
stem correlators, assuming that tails and arms are known in degrees $<d$,
and then reconstruct the arm $\tilde{T}(q)$ and tail $\delta t (q)$ 
(at $t=0$) up to degree $d$ from projections $J(0)_1$ and $J(0)_{\z}$ to 
$\K^{fake}_{+}$ and $\K^{\z}_{+}$ respectively. 
 
It is essential here that even when the head has degree $0$, it suffices to 
know the arms up to degree $<d$ (since at least 2 arms must be attached to 
the head). Also, when both the stem and the tail have degree 0, and there 
is only one leg attached, Proposition 3 recovers the information about 
the leg up to degree $d$ from that of the arm up to degree $d/m<d$. 

The previous procedure reconstructs $\tilde{T}$ (the arm at 
$t=0$), and hence the leg $T=\Psi^m(\tilde{T})$ in all degrees. Now, 
starting with any (non-zero) input $t$, one can first determine $\tilde{t}$ 
up to degree $d$ from stem correlators, assuming that tails are known in 
degrees $<d$, and then recover $\J(t)$ (and hence arms and tails) 
up to degree $d$. $\Box$

\medskip

Thus, to complete the proof of the theorem, it remains to show that the
tangent spaces from Proposition 2 coincide with the Lagrangian
spaces $\L^{\z}=\nabla_{\z}\T^{\z}$ described in the adelic formulation of
the theorem. This will be done in the next section.  

\section{Stems as stable maps to $X/\ZZ_m$}   

Let $\z \neq 1$ be a primitive $m$th root of 1, and let 
$\M_{0,n+2}^{X,d}(\z)$ denote a {\em stem space}. It is formed by stems 
of degree $d$, considered as quotient maps by the symmetry of order $m$ 
acting by $\z$ on the cotangent line at the horn of the covering curve.
It is a Kawasaki stratum in $\M_{0,mn+2}^{X,md}$. 

\medskip

{\tt Proposition 5.} {\em The stem space $\M_{0,n+2}^{X,d}(\z)$ is 
naturally identified with the moduli space 
$\M_{0,n+2}^{X/\ZZ_m,d}(g, 1,\dots, 1, g^{-1})$ 
of stable maps to the orbifold $X/\ZZ_m$.}

\medskip
 
{\tt Remark.} This Proposition refers to the GW-theory of orbifold target
spaces in the sense of Chen--Ruan \cite{ChR} and 
Abramovich--Graber--Vistoli \cite{AGV}. In particular, evaluations at marked 
points take values in the inertia orbifold, and notation of the moduli space
indicates the {\em sectors}, i.e. components of the inertia orbifold where
the evaluation maps land. In the case at hands the inertia orbifold is 
$X\times \ZZ_m$, and the string $(g, 1, \dots, 1, g^{-1})$, where $g$ is the 
generator of $\ZZ_m$, designates (in a way independent of $\z$) the sectors
of the marked points. 

\medskip

{\tt Proof.} The paper \cite{JK} by Jarvis--Kimura describes stable maps to 
the orbifold $point/\ZZ_m = B\ZZ_m$ in a way that can be easily 
adjusted to our case $X/\ZZ_m = X\times B\ZZ_m$. Namely, they are stable maps 
to $X$ equipped with a principal $\ZZ_m$-cover over the complement to the set 
of marked and nodal points, possibly ramified over these points in a way 
{\em balanced} at the nodes (i.e. such that the holonomies around the node
on the two branches of the curve are inverse to each other). The stem space
is obtained when two marked points are assigned holonomies $g^{\pm 1}$ of 
maximal order, and all other marked points are unramified. $\Box$

\medskip

Thus, introducing the simplifying notation $\M:=\M_{0,n+2}^{X,d}(\z)$, 
we identify stem correlators in the virtual Kawasaki formula \cite{To1} with 
integrals:
\begin{align*} \left[ \frac{}{}\right. \frac{\Phi}{1-q\z L^{1/m}} , & 
T(L),\dots,T(L), \delta t (L) \left. \frac{}{} \right]_{0,n+2}^{X,d} = \\  
\int_{[\M ]^{vir}} \td &(T_{\M})
\ch\left(\frac{\ev_1^*\Phi \ \ 
\ev_{n+2}^*\delta t (\z^{-1}L^{1/m}_{n+2})\ \ \prod_{i=2}^{n+1} \ev_i^*T(L_i) }
{(1-q\z L_1^{1/m})\ \Tr \left( \bigwedge^{\bullet} N^*_{\M}\right)} \right), 
\end{align*} 
Here $[\M]^{vir}$ is the virtual fundamental cycle of the moduli space
in GW-theory of $X/\ZZ_m$, $T_{\M}$ is the virtual tangent bundle to $\M$, and
$N_{\M}$ is the normal bundle to $\M$ considered as a Kawasaki stratum in the 
appropriate moduli space of stable maps to $X$ which are the $m$-fold covers of
the stems. In several steps, we will express stem correlators 
in terms of cohomological GW-theory of $X$. 

\medskip

Let $(\H, \Omega)$ be the symplectic loop space of cohomological GW-theory of 
$X$:
\[ \H = H ((z)),\ \ \Omega (f,g) = \Res_{z=0} (f(-z),g(z)) \ dz, \ \ 
(A,B)=\int_X AB.\]
Recall that the J-function of this theory is 
\[ \J^H_X(t):=-z+t(z)+\sum_a \phi_a \sum_{n,d} \frac{Q^d}{n!} 
\lan \frac{\phi^a}{-z-\psi}, t(\psi), \dots,t(\psi)\ran_{0,n+1}^{X,d}.\]
 
Consider now two generating functions for GW-invariants of $X/\ZZ_m$:
\[ \J_{X/\ZZ_m}^H(t):=-z+t(z)+\sum_a \phi_a \sum_{n,d}\frac{Q^d}{n!}
\lan \frac{\tilde{\phi^a}}{-z-\psi}, t(\psi),\dots,t(\psi)\ran_{0,n+1}
^{X/\ZZ_m, d}\]
and $\delta \J_{X/\ZZ_m}^H:=$ 
\[ \delta t(z)+\sum_a \phi_a\sum_{n,d}
\frac{Q^d}{n!}
\lan \frac{g\tilde{\phi^a}}{-z-\psi},t(\psi),\dots, t(\psi),g^{-1}
\delta t(\psi) \ran_{0,n+2}^{X/\ZZ_m, d}.\]
The factors $g^{\pm 1}$ indicate the appropriate sectors, while the 
inputs $t(\psi)$ are assumed to come from the identity sector. Here 
$\tilde{\phi^a}$ form a basis dual to $\phi_a$ with respect
to the Poincare pairing $A,B \mapsto (A,B)/m$. 

\medskip

{\tt Proposition 6.} {\em $\J_{X/\ZZ_m}^H(t) = \J_X^H(t)$, and 
\[ \delta \J_{X/\ZZ_m}^H =  \delta t(z)+\sum_a\phi_a\sum_{n,d}
\lan \frac{\phi^a}{-z-\psi},t(\psi),\dots,t(\psi),\delta t(\psi)
\ran_{0,n+2}^{X,d}.\]}

\medskip
  
{\tt Proof:} This follows directly from results of Jarvis--Kimura \cite{JK}.   
 
\medskip

{\tt Remark.} In fact the overruled Lagrangian cone of genus-0 GW-theory
on $X/\ZZ_m$ is the Cartesian product of $m$ copies of $\L^H_X$ corresponding
to {\em characters} of $\ZZ_m$. 
The J-function $\J^H_{X/\ZZ_m}$ represents points on 
the {\em diagonal} of this product. At such a point, the tangent space 
to the whole product decomposes into the direct sum of $m$ subspaces
(copies of $T_{J^H_X}\L^H_X$) according to {\em sectors} 
(i.e. elements of $\ZZ_m$). The vectors $\delta \J_{X/\ZZ_m}^H$ lie in one of 
these tangent subspaces, namely the one corresponding to the sector 
$g\in \ZZ_m$.
 
\medskip

The sum $T_{\M}\oplus N_{\M}$ is the restriction to $\M$ of the 
virtual tangent bundle to the moduli space of stable maps of degree 
$md$ with $mn+2$ marked points. According to \cite{Co}, in the Gr\"othendieck 
group $K^0(\M)$, this tangent bundle is represented by push-forward from the 
universal family $\tilde{\pi}: \tilde{\C} \to \M$:\footnote{In \cite{Co},
we have $T_X-1$ in place of $T_X$, but $\tilde{\pi}_*1=0$ in genus $0$.}   
\[ \tilde{\pi}_*\ev^*T_X+\tilde{\pi}_*(1-L^{-1})+\left(-\tilde{\pi}_*\tilde{i}_*
{\mathcal O}_{\tilde{Z}}\right)^{\vee},\]
where $L$ stands for the universal cotangent line at the ``current'' ($mn+3$-rd)
marked point of the universal curve, $\tilde{i}:\tilde{Z}\to \tilde{C}$ is the 
embedding of the nodal locus, and $\vee$ means dualization. 
This decomposes the virtual bundle into the sum of three parts, 
respectively responsible for: 
(i) deformations of maps to $X$ of a fixed complex curve, (ii) deformations 
of complex structure and/or configuration of marked points, and (iii) 
bifurcations of the curve's combinatorics through smoothing at the nodes. 

Part (i) is the {\em index bundle}
\[ \Ind (T_X) := \tilde{\pi}_*\tilde{\ev}^*(T_X). \]
Here we use the following notation: maps
$\pi: \C \to \M$ and $\ev:\C \to X/\ZZ_m$ form the universal stable map 
diagram, while $\tilde{\pi}: \tilde{\C} \to \M$ and 
$\tilde{\ev}: \tilde{\C} \to X$ are their $\ZZ_m$-equivariant lifts to the
family of ramified $\ZZ_m$-covers. 

We need to extract from the index bundle the eigenspace of the generator, 
$g$, of the group $\ZZ_m$, with the eigenvalue $\z^{-k}$.  For this, 
we begin with the $\ZZ_m$-module $\CC$ where $g$ acts by $\z^{k}$, denote
$\CC_{\z^{k}}$ the corresponding line bundle over $B\ZZ_m$, and take 
$\left( \Ind (T_X) \otimes \CC_{\z^{k}} \right)^{\ZZ_m}$. 
This (trivial) result can be expressed in terms of orbifold GW-theory 
of $X/\ZZ_m=X\times B\ZZ_m$ as $\pi_*\ev^*\left(T_X \otimes \CC_{\z^{k}}
\right)$. Namely, as we mentioned in Section 4, the K-theoretic push-forward 
operation on global quotients considered as orbifolds automatically 
extracts the invariant part of sheaf cohomology. Thus, 
\[ \Tr (\Ind (T_X)) = \oplus_{k=0}^{m-1} \z^{-k} \pi_*\ev^*\left(T_X\otimes 
\CC_{\z^{k}}\right).\]

Recall that an invertible characteristic class of complex vector bundles is
determined by an invertible formal series in one variable, the 1st Chern
class $x=c_1(l)$ of the universal line bundle. Alongside the usual Todd class
$\td$, we introduce {\em moving} Todd classes (aka equivariant K-theoretic 
inverse Euler classes), one for each $\lambda \neq 1$: 
\[ \td (l) = \frac{x}{1-e^{-x}}, \ \ \ \td_{\lambda}(l) =
\frac{1}{1-\lambda e^{-x}}.\]
The contribution of $\Tr (\Ind (T_X))$ into our integral
over $\M$ reads:
\begin{equation} \label{twisting_factor} \tag{*}
\td \left(\pi_*\ev^*(T_X) \right) \ \prod_{k=1}^{m-1} \td_{\z^k} \left(
\pi_*\ev^*(T_X\otimes \CC_{\z^k})\right) .\end{equation}
Introduce $\J^{tw}_{X/\ZZ_m}$ and $\delta \J^{tw}_{X/\ZZ_m}$ as {\em twisted} 
counterparts of $\J^H_{X/\ZZ_m}$ and $\delta \J^H_{X/\ZZ_m}$. Namely,     
following \cite{CGi}, one defines GW-invariants {\em twisted} by a chosen
bundle, $E$, over the target space, and a chosen multiplicative 
characteristic class, $S$, by systematically replacing virtual fundamental 
cycles of moduli spaces of stable maps with their cap-products 
(such as $[\M]^{vir} \cap S(\Ind (E))$ in our case) with the chosen
characteristic class of the corresponding index bundle. 

\medskip

{\tt Proposition 7.} {\em Denote by $\square$ and $\square_{\z}$
the Euler--Maclaurin asymptotics of the infinite products 
\begin{align*} 
 \square &\sim
\prod_{\text{Chern roots $x$ of $T_X$}}\ \  \prod_{r=1}^{\infty} 
\frac{x-rz}{1-e^{-mx+mrz}} \\
\square_{\z} &\sim  
\prod_{\text{Chern roots $x$ of $T_X$}} \ \ \prod_{r=1}^{\infty} 
\frac{x-rz}{1-\z^{-r} e^{-x+rz/m}}. \end{align*}
Then $\J^{tw}_{X/\ZZ_m}$ lies in the overruled Lagrangian cone $\square 
\L^H_X$, and $\delta\J^{tw}_{X/\ZZ_m}$ lies in the transformed tangent space 
$\square_{\z} \T_{\square^{-1}\J^{tw}_{X/\ZZ_m}} \L^H_X$.} 

\medskip

{\tt Proof.} The Quantum
Riemann--Roch Theorem of \cite{CGi}, which expresses twisted GW-invariants
in terms of {\em un}twisted ones, was generalized to the case of orbifold
target spaces by Hsian-Hua Tseng \cite{Ts}. The proposition is obtained by 
direct applications of the Quantum RR Theorem of \cite{Ts} to each of the 
twisting data $E=T_X \otimes \CC_{\z^k}$, $S=\td_{\z^k}$. For $k=0$, the 
Euler--Maclaurin asymptotics (for both
$\J^{tw}$ and $\delta \J^{tw}$) come from the product
\[ \prod_{r=1}^{\infty} \frac{x-rz}{1-e^{-x+rz}}, \]
and for $k\neq 0$, from
\[ \prod_{r=1}^{\infty} \frac{1}{1-\z^{-k} e^{-x+rz}}\]
for $\J^{tw}$, and
\[ \prod_{r=1}^{\infty} \frac{1}{1-\z^{-k} e^{-x+rz+kz/m}} \]
for $\delta \J^{tw}$. 
Multiplying out the products over $k=0,\dots,m-1$, and using 
$\prod_{k=0}^{m-1} (1-\z^{-k}u)=1-u^m$ and $\z^m=1$
to simplify, we obtain the required results. $\Box$

\medskip

Part (ii) of the bundle $T_{\M}\oplus N_{\M}$ comes from deformations of the
complex structure and marked points. It can be described as the K-theoretic 
push-forward $\tilde{\pi}_*(1-L^{-1})$ 
along the universal curve $\tilde{\pi}: \tilde{\C} \to \M$ 
(think of $H^1(\Sigma, T_{\Sigma})$). To express the trace $\Tr$ of it,
one need to consider push-forwards of $L^{-1}\otimes \CC_{\z^k}$ and 
appropriately twisted GW-invariants of the orbifold $X/\ZZ_m$. More precisely,
we need the twisting classes to have the form:
\[ \td \left(\pi_* (1-L^{-1})\right) \ \prod_{k=1}^{m-1}\td_{\z^k}\left(\pi_*
[(1-L^{-1})\ev^*(\CC_{\z^k})]\right).\]  
The general problem of computing GW-invariants of orbifolds twisted 
by characteristic classes of the form
\[ \prod_{\a} S_{\a} \left(\pi_*[(L^{-1}-1)\ev^*(E_{\a})]\right)\]  
is solved in \cite{To2} (see also Chapter 2 of thesis \cite{To}). 
The answer is described as the change of the dilaton shift.\footnote{
Generalizing the case of manifold target spaces discussed in \cite {Co}.}
Namely, if $-z=c_1(L^{-1})$, and $S_{\a}$ denote the twisting multiplicative 
characteristic class, then the dilaton shift changes from $-z$ to 
$-z \prod_{\a} S_{\a}(L^{-1} E_{\a})$.   
In our situation, $\a=0,\dots ,m-1$, $S_0=\td^{-1}$, $S_k=\td_{\z}^{-1}$ for 
$k\neq 0$, and $E_k=\CC_{\z^k}$. Respectively, the new dilaton shift is
\[ -z \frac{(1-e^z)}{(-z)} \prod_{k=1}^{m-1}(1-\z^k e^{z})=1-e^{mz}.\]
Thus, the {\em dilaton shift changes from $-z$ to $1-e^{mz}$}. 

\medskip

Parts (i) and (ii) together form the part of the virtual tangent bundle to
$\M_{0,mn+2}^{md,X}$ (albeit restricted to $\M$) {\em logarithmic} with 
respect to the nodal divisor. What remains is part (iii), supported on the
nodal divisor, which consists of one-dimensional summands (one per node), 
the {\em smoothing mode} of the glued curve at the node. Contributions of part
(iii) into the ratio $\td (T_{\M})/ \ch (\Tr \bigwedge^{\bullet}(N^*_{\M})$ 
in the virtual Kawasaki formula is described in terms of yet another kind of
twisted GW-invariants of the orbifold $X/\ZZ_m$, where the twisting classes
are supported at the nodal locus.  The effect of such twisting on GW-invariants
can be found by a recursive procedure based on {\em un}gluing the curves at 
the nodes. As it is seen in \cite{Co}, this does not change the overruled 
Lagrangian cones, but affects generating functions through a {\em change of 
polarization}. Referring to \cite{To2} (or \cite{To}) for the generalization 
to orbifold target spaces needed here, we state the results.
     
Let $\M$ denote a moduli space of stable maps to the orbifold $X/\ZZ_m$, and
$\pi: \C\to \M$ the projection of the universal family of such stable maps.
Let $Z=\cup_{h\in \ZZ_m} Z_h$ be the decomposition of the nodal stratum
$Z \subset \C$ into the disjoint union according to the ramification type of 
the node, and $i:Z_h\to \C$ denote the embedding.
Let $S_{h,a}$ be invertible multiplicative 
characteristic classes, and $E_{h,\a}$ arbitrary orbibundles 
over $X/\ZZ_m$, where $h\in \ZZ_m$, $\a=1,\dots, K_h$.
The twisting in question is obtained by systematically 
including into the integrands of GW-theory of $X/\ZZ_m$ the factors
\[ \prod_{h\in \ZZ_m} \prod_{\a=1}^{K_h} S_{h,\a}\left(\pi_*[i_*{\mathcal O}_{Z_h} 
\otimes \ev^*E_{h,\a}]\right).\]
According to the results of \cite{To2} (Theorem 1.10.3 in \cite{To}), the 
effect of such twisting is completely accounted by a change of polarization 
in the symplectic loop space of GW-theory of $X/\ZZ_m$, described separately
for each sector. Namely, for the sector corresponding to $h\in \ZZ_m$, define
a power series $u_h(z)=z+a_2z^2+a_3z^3+\cdots$ by
\[ \frac{z}{u_h(z)} = \prod_{\a=1}^{K_h} S_{h.a}^{-1} \left(E_{h,a}\otimes 
L\right),\ \ \text{where $c_1(L):=z$.}\]
Define the Laurent series $v_{h,k}, k=0,1,2,\dots$, by
\[ \frac{1}{u_h(-\psi-z)}=\sum_{k\geq 0} \left(u_h(\psi)\right)^k v_{h,k}(z),\]
which is the expansion of the L.H.S. in the region $|\psi|<|z|$. Then, as one
can check,  
$\phi_{a} z^k, \phi^{a}v_{h,k}(z)$, $a=1,\dots, \dim H$, $k=0,1,2,\dots$, 
form a (topological) Darboux basis in the sector $h$ of the symplectic 
loop spaces of the GW-theory of $X/\ZZ_m$. The genus-0 descendant potential 
of the twisted theory is expressed from that of untwisted one by taking the 
overruled Lagrangian cone of the untwisted theory for the graph of 
differential of a function in the Lagrangian polarization associated with 
this basis. Note that the positive space polarization, which is
spanned by $\{ \phi_{a} z^k \}$, stays the same as in the untwisted theory, 
while the negative space, which is spanned by $\{ \phi^{a} v_{h,k}(z)\}$,
differs from that of untwisted theory, which is spanned by 
$\{ \phi^{a} z^{-1-k}\}$.  

\medskip

{\tt Remarks.} (1) The standard polarization $\H_{\pm}$ of the symplectic loop
space of quantum cohomology theory of a manifold is obtained by the same 
formalism:
\[ \frac{1}{-\psi-z}=\sum_{k\geq 0} \frac{\psi^k}
{(-z)^{k+1}},\] 
and $\H_{-}$ is spanned by $\phi^a (-z)^{-1-k},\ k=0,1,2,\dots$.

(2) As it was mentioned in Section 5, in fake K-theory one obtains a Darboux 
basis from $z/u(z)=\td (L)$, and respectively the expansion:
\[ \frac{1}{1-e^{\psi+z}} = 
\sum_{k\geq 0} (e^{\psi}-1)^k\frac{e^{kz}}{(1-e^{z})^{k+1}}
.\]
Consequently, $\K_{+}^{\fake}$ and $\K_{-}^{\fake}$ are spanned respectively by 
$\Phi_a (q-1)^k$ and $\Phi^a q^k/(1-q)^{k+1}$, $a=1,\dots,\dim K$, 
$k=0,1,2,\dots$.
    
\medskip

In stem theory, there are two types of nodes (Figure 2). 
When a stem acquires an unramified node (as shown in the top picture), 
the covering curve carries a $\ZZ_m$-symmetric $m$-tuple of nodes. 
The smoothing bundle has dimension $m$ and carries a regular representation 
of $\ZZ_m$. When a stem degenerates into a chain of two components glued 
at a balanced ramification point of order $m$ (the bottom picture), the 
smoothing mode is one-dimensional and carries the trivial representation of
$\ZZ_m$.  Contributions of these smoothing modes into the ratio 
$\td (T_{\M})/ \ch (\Tr \bigwedge^{\bullet}(N^*_{\M})$ is accounted by 
the following twisting factors in the integrals over $\M$, considered 
as orbifold-theoretic GW-invariants:
\[ \td \left(-\pi_*i_*{\mathcal O}_{Z_g}\right)^{\vee} 
\td \left(-\pi_*i_*{\mathcal O}_{Z_1}\right)^{\vee} \prod_{k=1}^{m-1} 
\td_{\z^k} \left(-\pi_*(\ev^*\CC_{\z^k} \otimes i_*{\mathcal O}_{Z_1})\right)^{\vee},
\] 
where $Z_1$ stands for the unramified nodal locus, and $Z_g$ for the ramified 
one. This twisting results in the change of polarizations. In the $g$-ramified
sector, the new polarization is determined from the expansion of
\[ \frac{1}{1-e^{(psi+z)/m}} = \frac{1}{1-q^{1/m}L^{1/m}}.\]
Here the factor $1/m$ occurs because what was denoted $L$ in the GW-theory
of $X/\ZZ_m$ is the universal cotangent line at the ramification point to the
quotient curve, which is $L^{1/m}$ in our earlier notations of stem spaces
(where $L$ stands for the universal cotangent line to the covering curve). 
In the unramified sector, the new polarization is found from
\[ \frac{z}{u(z)}=\td (L) \prod_{k=1}^{m-1}\td_{\z^k}(\CC_{\z^k}\otimes L) = 
\frac{z}{1-e^{-z}}\prod_{k=1}^{m-1}\frac{1}{1-\z^ke^{-z}} = \frac{z}{1-e^{-mz}},\]
and consequently the expansion of
\[ \frac{1}{1-e^{m\psi+mz}}=\frac{1}{1-q^mL^m}.\]
We conclude that the {\em negative space of polarizations in the ramified and
unramified sectors are spanned respectively by}
\[ \Phi^a q^{k/m}/(1-q^{1/m})^{k+1} \ \ \text{and}\ \  
\Phi^a q^{mk}/(1-q^m)^{k+1} =\Phi^a\ \Psi^m\left( q^k/(1-q)^{k+1}\right).\]
 
\medskip

{\tt Remark.} The occurrence of Adams' operation $\Psi^m$ here is not
surprising. The smoothing modes at $m$ cyclically permuted copes on
an unramified node of the stem curve form an $m$-dimensional space carrying 
the regular representation of $\ZZ_m$. The trace $\Tr$ of the bundle formed
by these modes is, according to Lemma of the previous section, 
$\Psi^m (L_{-}\otimes L_{+})$ (in notation of Figure 2, the top picture). 
 
\medskip

\begin{figure}[htb]
\begin{center}
\epsfig{file=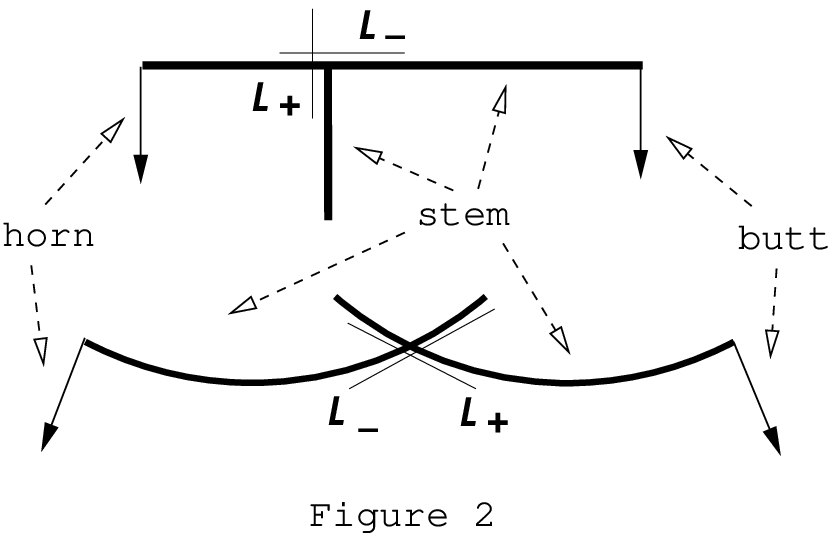}
\end{center}
\end{figure}    

It remains to apply the above results to those generating 
functions of stem theory which occur in the virtual Kawasaki formula.
Introduce a generating function, 
$\delta \J^{st}_{X/\ZZ_m}$, of stem theory as the image under the Chern 
character map $\ch: \K^{\fake} \to \H$ of 
\[ \delta{t}(q^{1/m})+
\sum_{a, n,d}\Phi_a \frac{Q^{d}}{n!}\left[ \frac{\Phi^a}
{1-q^{1/m}L^{1/m}}, T(L),\dots,T(L),\delta t(L^{1/m})
\right]_{0,n+2}^{X,d}. \]
Replacing $q^{1/m}$ with $\z q$ and $Q^d$ with $Q^{md}$, we would 
obtain the sum of correlators of stem theory as they appeared in Section 7. 
On the other hand, the interpretation of stem correlators as GW-invariants
of $X/\ZZ_m$ twisted in three different ways (corresponding to parts 
(i), (ii), (iii) of the tangent space), and the previous results on the effects
of these twistings, provide the following description of
$\delta \J^{st}{X/\ZZ_m}$ in terms of GW-theory of $X$. 

\medskip

{\tt Proposition 8.} {\em $\delta \J^{st}_{X/\ZZ_m}
(\delta t, T)$  lies in 
$\square_{\z} \square^{-1} \T_{\J^{tw}_{X/\ZZ_m}} \square \L^H_X$,
where
the input $T$ is related to the application point 
$\J^{tw}_{X/\ZZ_m}$ by the projection $\left[ \cdots \right]_{+}$
along the negative space of the polarization of the unramified 
sector:
\[ \ch \left( 1-q^m + T(q) \right) = 
\left[ \J^{tw}_{X/\ZZ_m} \right]_{+} \ .\]}

%\medskip

{\tt Proof.} According to Proposition 7, $\J^{tw}_{X/\ZZ_m}$ lies in the cone
$\square \L^H$, and $\delta J^{tw}_{X/\ZZ_m}$ lies 
in the space $\square_{\z} \square^{-1} \T$, where $\T$ is the tangent
space to $\square \L^H$ at the point $\J^{tw}_{X/\ZZ_m}$. It follows from the
previous discussion that $\delta \J^{st}_{X/\ZZ_m}$, being obtained from 
$\delta \J^{tw}_{X/\ZZ_m}$ by changing dilaton shift and polarizations only, 
lies in the same space. Changing the content of the horn 
in the definition of $\delta J^{tw}_{X/\ZZ_m}$ from $\tilde{\phi}^a/(-z-\psi)=
\phi^a/(-z/m-\psi/m)$ to $\phi^a/(1-e^{z/m+\psi/m})$ is equivalent to 
applying to the same space the polarization associated with the $g$-ramified
sector. However, the new dilaton shift and polarization in the unramified
sector both affect the way the input $T$ 
of $\delta \J^{st}_{X/\ZZ_m}$ is computed in terms of $\J^{tw}_{X/\ZZ_m}$. 
Namely, $\ch(T)=[\J^{tw}_{X/\ZZ_m}]_{+} - (1-e^{mz})$. $\Box$
          
\medskip

To put the next proposition into a context, let us recall that 
the stem correlators of Section 7, in order
to represent the expansion $\J(t)_{\z}$ of the true K-theoretic J-function 
$\J$ at $q=\z^{-1}$, need to
be computed at a specific input $T$, the {\em leg}, which is characterized
in a rather complex way. Namely, the expansion $\J(0)_1$ of the value of 
$\J$ at the input $t=0$ lies in the cone $\L^{\fake}$ of fake quantum K-theory
of $X$ (Proposition 1). The contribution $\tilde{T}$, i.e. the {\em arm}
corresponding to $t=0$, is obtained as the input point of $\J(0)_1$, i.e. 
by applying the projection $(\dots)_{+}$ along the negative space 
of polarization (described in Example 2) and the dilaton shift of 
fake quantum K-theory:
\[ 1-q+\tilde{T}(q) = \left(\J(0)_1\right)_{+}.\]
According to Proposition 3, $T=\Psi^m(\tilde{T})$ (where Adams' operation acts
also on $q$ and $Q$).

On the other hand, Proposition 7 locates the stem generating function 
$\delta\J^{st}_{X/\ZZ_m}$ in terms of the tangent space 
$\T_{\J^{tw}_{X/\ZZ_m}} \square \L^H_X$. Furthermore, according to the Quantum HRR 
Theorem stated in Section 5, $\L^H_X=\Delta^{-1} \ch\left(\L^{\fake}\right)$, 
and hence the cone and its tangent space in question are the images under 
$\square \Delta^{-1}\ch$ of $\L^{\fake}$ and of a certain tangent space to it. 
The following proposition implies that when the input $T$ is the leg, the 
the requisite tangent space is exactly $T_{\J(0)_1}\L^{\fake}$.
      
\medskip
 
{\tt Proposition 9.} {\em $\ch^{-1} (\square \L^H_X) = \Psi^m (\L^{\fake})$, 
where the Adams operation $\Psi^m: \K^{\fake}\to \K^{\fake}$ acts on 
$q$ by $\Psi^m (q) := q^m$.}
%{\em We have:
%\[ \J^{tw}_{X/\ZZ_m} = m^{\frac{1}{2}\dim_{\CC} X-1}\  
%\ch \left( \Psi^m (\J_X^{\fake} ) \right),\]
%where the Adams operation $\Psi^m$ acts on $q$ by $\Psi^m (q) := q^m$.}
  
\medskip

{\tt Proof.} From the QHRR theorem of Section 5 and Proposition 7, we have: 
\[ \Delta^{-1}\ch (\J_X^{\fake}) = \J^H_X = 
\square^{-1} \J^{tw}_{X/\ZZ_m}.\]
We intentionally neglect to specify the arguments, since they are determined by
the argument, $t$, of $\J^H_X$, by polarizations, and by the 
transformations $\Delta$ and $\square$ themselves. The Adams operation 
$\Psi^m$ acts on cohomology classes via the Chern isomorphism:
\[ \ch \left( \Psi^m (\ch^{-1} a) \right) = m^{\deg (a)/2} a.\] 
The J-function $\J^H_X$ has degree $2$ with respect to the grading, defined by
$\deg z=2$, $\deg Q^d = 2 \int_d c_1(T_X)$, and the usual grading in cohomology. 
Therefore, writing $\J^H_X/(-z) = \sum_d J_d Q^d$, we find 
\[ m^{-1}\Psi^m (\J^H_X) = \sum_d m^{-\deg Q^d} J_d Q^d = 
e^{-(\log m)\ c_1(T_X)/z} \J_X^H. \] 
The last equality is an instance of the genus-0 {\em divisor equation} (see 
\cite{CGi}). Thus, Proposition 9 would follow from the identity
\[ \square =m^{\frac{1}{2}\dim_{\CC} X} \Psi^m(\Delta)\ e^{-(\log m)\ 
c_1(T_X)/z}.\]  
To establish it, note that both $\Delta$ and $\square$ are the 
Euler--Maclaurin asymptotics of infinite products 
\[ \prod_{\text{Chern roots $x$ of $T_X$}}\ \  \prod_{r=1}^{\infty} S(x-rz),\]
where $S$ is respectively \[ \frac{x}{1-e^{-x}}\ \ \ \ \text{and}\ \ \ \  
 \frac{x}{1-e^{-mx}} = m^{-1}\ \Psi^m \left( \frac{x}{1-e^{-x}} \right).\] 
The factor $m^{-1}$ contributes into the asymptotics in the form
\[ \prod_{\text{Chern roots $x$ of $T_X$}}e^{-(\log m)\ x/z}m^{1/2} 
= e^{-(\log m)\ c_1(T_X)/z} 
m^{\frac{1}{2}\dim_{\CC} X}. \hspace{15mm} \Box \]

\medskip

{\tt Corollary.} {\em  
$\delta \J^{st}_{X/\ZZ_m} (\delta t, T)$ lies in the space 
$\square_{\z}\Delta^{-1}  \T_{\J^{\fake}(\tilde{T})} \L^{\fake}$, 
where $T=\Psi^m(\tilde{T})$.} 

{\tt Proof.} According to Proposition 8, 
$\delta \J^{st}_{X/\ZZ_m}$ lies in the space \linebreak
$\square_{\z}\square^{-1}\T_{\J^{tw}_{X/\ZZ_m}}\square \L^H$, where the 
input $T$ of $\delta \J^{st}_{X/\ZZ_m}$ is determined by 
$T=[\J^{tw}_{X/\ZZ_m}]_{+}-(1-q^m)$. By Proposition 9, 
$\J^{tw}_{X/\ZZ_m}=\Psi^m (\J^{\fake})$, and the input of $\J^{\fake}$
is determined as $\tilde{T}=\left(\J^{\fake}\right)_{+}-(1-q)$. Here 
$\left(\cdots\right)_{+}$ refers to the projection to $\K^{\fake}_{+}$
along $\K^{\fake}_{-}$, i.e. the polarization described in Example 2,
while the projection $[\cdots]_{+}$ refers to the polarization in the unramified
sector. The latter polarization is obtained from the former 
by the Adams operation: $\Psi^m: \K^{\fake} \to \K^{\fake}$, 
and the relation between dilaton shifts is the same:
$1-q^m=\Psi^m (1-q)$. Therefore $T=\Psi^m(\tilde{T})$, and the tangent space
$\T_{\J^{tw}_{X/\ZZ_m}}\square \L^H$ can be described as 
$\square \Delta^{-1} \T_{\J^{\fake}(\tilde{T})} \L^{\fake}$.
$\Box$

\medskip

We note that 
\[ \square_{\z} \Delta^{-1} \sim \prod_{\text{Chern roots $x$ of $T_X$}}
\prod_{r=1}^{\infty} \frac{1-q^r e^{-x}}{1-\z^{-r}q^{r/m} e^{-x}}.\]
One obtain $\nabla_{\z}$ by replacing in this formula $q^{1/m}$ with $q\z$ and
computing the Euler--Maclaurin asymptotics of the result as $q\z \to 1$. 

%Recall that $\J(0)_1$, the expansion near $q=1$ of the value $\J(t)$
%at the input $t$ of the big J-function of true quantum K-theory of $X$,
%lies in $\L^{\fake}$ (Proposition 1). According to Propositions 2 and 3, 
%$\J(t)_{\z}$ (the expansion of $\J(t)$ near $q=\z^{-1}$) is given by the 
%generating function of stem theory with the contribution 
%$T(L)$ of the leg obtained from the contribution $\tilde{T}(L)$ of the arm 
%without marked points (i.e. from the input of $\J(0)_1$ considered as a 
%value of $\J^{\fake}$) by the operation $\Psi^m$ and the change 
%$Q^d\mapsto Q^{md}$. 

According to Propositions 2, the expansion $\J(t)_{\z}$ near $q=\z^{-1}$ 
of the true K-theoretic J-function is expressed in terms of correlators 
of stem theory (as they appeared in Section 7), computed at the input $T$ 
equal to the leg contribution. More precisely, $J(t)_{\z}$ is obtained
from $\delta \J^{st}_{X/\ZZ_m}$, defined as:  \[ \delta{t}(q^{1/m})+
\sum_{a, n,d}\Phi_a \frac{Q^{d}}{n!}\left[ \frac{\Phi^a}
{1-q^{1/m}L^{1/m}}, T(L),\dots,T(L),\delta t(L^{1/m})
\right]_{0,n+2}^{X,d}, \]
by changing $q^{1/m} \to q\z$ (including such change in $\delta t$) and 
$Q^d\to Q^{md}$ (excluding such a change in $\delta t$).  
   
Combining these facts with Corollary, we conclude that  
$\nabla_{\z}^{-1}\J(t)_{\z}$, after the change $q\mapsto q/\z$, falls into 
the subspace $\T$ of $\K^{\fake}$ which is obtained
from the tangent space $\T_{\J(0)_1}\L^{\fake}$ by the changes 
$q^{1/m}\mapsto q$ and $Q^d\mapsto Q^{md}$. 

\medskip

This completes the proof of the Hirzebruch--Riemann--Roch Theorem
in true quantum K-theory. 
 
\section{Floer's $S^1$-equivariant K-theory, and $\D_q$-modules}

In this section, we show that tangent spaces to the overruled Lagrangian
cone $\L$ of quantum K-theory on $X$ carry a natural structure of modules
over a certain algebra $\D_q$ of finite-difference operators with respect 
to Novikov's variables. This structure, although manifest in some examples 
(see \cite{GiL}) and predictable on heuristic grounds of $S^1$-equivariant 
Floer theory \cite{GiZ, GiH}, has been missing so far in the realm of 
K-theoretic GW-invariants. We first recall the heuristics, and then derive the
$\D_q$-invariance of the tangent spaces to $\L$ from the divisor equation 
in quantum cohomology theory and our HRR Theorem in quantum K-theory.

\medskip

Let $X$ be a compact symplectic (or K\"ahler) target space, which for 
simplicity is assumed simply-connected in this preliminary discussion, 
so that $\pi_2(X)=H_2(X)$. Let $k=\operatorname{rk} H_2(X)$, 
let $d=(d_1,\dots,d_k)$ be integer coordinates on 
$H_2(X,\QQ)$, and $\omega_1, \dots, \omega_k$ be closed 2-forms on $X$
with integer periods, representing the corresponding basis of $H^2(X,\RR)$.

On the space $L_0X$ of contractible parametrized loops $S^1\to X$, as well 
as on its universal cover $\tilde{L_0X}$, one defines closed 2-forms 
$\Omega_a$, that to two vector fields $\xi$ and
$\eta$ along a given loop associates the value
\[ \Omega_a(\xi,\eta) := \oint \omega_a(\xi(t),\eta(t))\ dt.\]  
A point $\gamma \in \tilde{L_0X}$ is a loop in $X$ together with a homotopy
type of a disk $u: D^2\to X$ attached to it. One defines the 
{\em action functionals} $H_a: \tilde{L_0X}\to \RR$ by evaluating 
the 2-forms $\omega_a$ on such disks:
\[ H_a(\gamma) := \int_{D^2} u^*\omega_a .\]
  
Consider the action of $S^1$ on $\tilde{L_0X}$, defined by the rotation 
of loops, and let $V$ denote the velocity vector field of this action.
It is well-known that $V$ is
$\Omega_a$-hamiltonian with the Hamilton function $H_a$, i.e.:
\[ i_{V} \Omega_a + dH_a = 0,\ \ \ a=1,\dots, k.\]

Denote by $z$ the generator of the coefficient ring $H^*(BS^1)$ of 
$S^1$-equivariant cohomology theory. The $S^1$-equivariant De Rham complex
(of $\tilde{L_0X}$ in our case) consists of $S^1$-invariant differential 
forms with coefficients in $\RR [z]$, and is equipped with the differential
$D:= d+zi_V$. Then 
\[ p_a:=\Omega_a + z H_a, \ \ \ \ a=1,\dots, k,\]
are degree-2 $S^1$-equivariantly closed elements of the complex:
%\linebreak 
$Dp_a=0$. This is a standard fact that usually accompanies the formula of
Duistermaat--Heckman.     

Furthermore, the lattice $\pi_2(X)$ acts by deck transformations on the 
universal covering $\tilde{L_0X}\to L_0X$. Namely, an element $d\in \pi_2(X)$
acts on $\gamma \in \tilde{L_0X}$ by replacing the homotopy type $[u]$ of the
disk with $[u]+d$. We denote by $Q^d=Q_1^{d_1}\cdots Q_k^{d_k}$ the operation 
of pulling-back differential forms by this deck transformation. It is an 
observation from \cite{GiZ,GiH} that the operations $Q_a$ and the operations
of exterior multiplication by $p_a$ do not commute:
\[ p_a Q_b - Q_b p_a = -z Q_a \delta_{ab}.\]
These are commutation relations between generators of the algebra of 
differential operators on the k-dimensional torus:
\[ \left[ -z \p_{\tau_a} , e^{\tau_b} \right] = -z e^{\tau_a} \delta_{ab}.\] 
Likewise, if $P_a$ denotes the $S^1$-equivariant line bundle on $\tilde{L_0X}$
whose Chern character is $e^{-p_a}$, then tensoring vector bundles by $P_a$ 
and pulling back vector bundles by $Q_a$ do not commute:
\[ P_a Q_b = q Q_a P_b \delta_{ab}.\]
These are commutation relations in the algebra of finite-difference
operators, generated by multiplications and translations:
\[ Q_a\mapsto e^{\tau_a}, \ \ P_a\mapsto e^{z\p_{\tau_a}}=q^{\p_{\tau_a}}, 
 \ \ \text{where}\ \ q=e^z.\] 
Thinking of these operations acting on $S^1$-equivariant Floer theory of the 
loop space, one arrives at the conclusion that $S^1$-equivariant Floer 
cohomology (K-theory) should carry the structure of a module over
the algebra of differential (respectively finite-difference) operators.
Here is how this heuristic prediction materializes in GW-theory.  

\medskip   
  
{\tt Proposition 10.} {\em Let $\D$ denote the algebra of differential 
operators generated by $p_a, a=1,\dots,k$, and $Q^d$, with $d$ lying in the 
Mori cone of $X$. Define a representation of $\D$ on the symplectic loop 
space $\H=H^*(X,\CC[[Q]]) \otimes \CC ((z))$ using the operators 
$p_a-zQ_a\p_{Q_a}$ (where $p_a$ acts by multiplication in the classical 
cohomology algebra of $X$) and $Q^d$ (acting by multiplication in the Novikov 
ring). Then tangent spaces to the overruled Lagrangian cone $\L^H\subset \H$
of cohomological GW-theory on $X$ are $\D$-invariant.}

\medskip

{\tt Proof.} Invariance of with respect to multiplication by $Q^d$ is 
tautological since the Novikov ring $\QQ [[Q]]$ (which contains 
the semigroup algebra of the Mori cone: we assume that $d_a=\int_d p_a \geq 0$
for all $a$ and all $d$ in the Mori cone) 
is considered as the ground ring of scalars. To prove invariance
with respect to operators $p_a-zQ_a\p_{Q_a}$, recall from \cite{GiF} 
that tangent spaces to $\L^H$ have the form $S^{-1}_{\tau}\H_{+}$, 
where $H\ni \tau \mapsto S_{\tau}(z)$ is a matrix power series in $1/z$ whose 
matrix entries are the following cohomological GW-invariants:   
\[  S_a^b=\delta_a^b+\sum_{l,d}\frac{Q^d}{l!}\sum_{\mu} 
\lan \phi^a, \tau, \dots, \tau, \frac{\phi_b}{z-\psi}\ran_{0,n+2}^{X,d}.\]
The matrix $S_{\tau}$ lies in the {\em twisted loop group}, 
i.e. $S^{-1}_{\tau}(z)=S^*_{\tau}(-z)$. Let $\p_{\tau_a}$ denote the
differentiation in $\tau$ in the direction of the degree-2 cohomology class 
$p_a$. According to the {\em divisor equation},
\[ z Q_a\p_{Q_a} S_{\tau}(z)+S_{\tau}(z)p_a = 
z\p_{\tau_a} S_{\tau}(z). \]
In fact $z\p_{\tau_a} S = p_a \bullet S$,  
where $\bullet$ stands for quantum cup-product. (This follows from
the property of $\L^H$ to be overruled.) Transposing, we get:  
\[ (p_a-z Q_a\p_{Q_a}) S^{-1}_{\tau}(z)
=-z\p_{\tau_a} S^{-1}_{\tau}(z) =S^{-1}_{\tau}(z) (p_a\bullet). \]
Also, if $\tau=\sum_{\mu} \tau_{\mu} \phi_{\mu} \in H$, then for any $\mu$
\[ z\p_{\tau_{\mu}} S_{\tau}(z)= (\phi_{\mu}\bullet) S_{\tau}(z),\ \ 
\text{and hence}\ \ 
-z \p_{\tau_{\mu}} S^{-1}_{\tau}(z)= S^{-1}_{\tau} (\phi_{\mu}\bullet).\]
Thus, if $\tau=\sum \tau_{\mu}(Q)\phi_{\mu}$ and $h\in \H_{+}$, so that 
$f(z,Q)=S^{-1}_{\tau}(z) h(z,Q)$ lies in $\T_{\tau}$, then
\[ (p_a-zQ_a\p_{Q_a}) f = S^{-1}_{\tau}(z) \left[ 
(p_a\bullet) - z Q_a\p_{Q_a} + z\sum Q_a\p_{Q_a}\tau_{\mu}\ (\phi_{\mu}\bullet) \right] h.\]
Since $\H_{+}$ is invariant under the operator in brackets, the result follows.
$\Box$   
    
\medskip

{\tt Remarks.} (1) Each ruling space $z\T_{\tau}$,
and therefore the whole cone $\L^H$, is $\D$-invariant, too.

(2) Symbols of differential operators annihilating all columns of $S$ 
provide relations between operators $p_a\bullet$ in the quantum cohomology 
algebra of $X$ (see \cite{GiG}). 

\medskip

{\tt Corollary 1.} {\em Tangent and ruling spaces of $\L^{fake}$ are 
$\D$-invariant.}

\medskip

{\tt Proof.} In the QHRR formula $\ch (\L^{\fake}) = \Delta \L^H$ of Section 5,
the operator $\Delta$ commutes with $\D$, since it does not involve Novikov's
variables, and since the operators (which do occur in $\Delta$) 
of multiplication in the classical cohomology algebra of $X$ commute with 
$p_a$. $\Box$     

\medskip

{\tt Lemma.} {\em The subspace $\T\subset \K^{\fake}$ obtained from
$\T_{\J(0)_1}\L^{\fake}$ by the change $z\mapsto mz, Q\mapsto Q^m$, is
$\D$-invariant.}

\medskip

{\tt Proof.} The tangent space in question is 
$\Delta (z) S^{-1}_{\tau^{(0)}(Q)}(z,Q) \H_{+}$ for some $\tau^{(0)} =
\sum_{\mu} \tau^{0}_{\mu} \phi_{\mu} \in H$. 
(Recall that $\H_{+}=H[[z]]$, and $H=H^*(X,\C[[Q]])$.) The space $\T$ 
is therefore $\Delta (mz) S^{-1}_{\tau^{(0)}(Q^m)}(mz,Q^m) \H_{+}$,
where $\H_{+}$ is $\D$-invariant, and $\Delta$ commutes with $\D$. 
Since $zQ_a\p_{Q_a}=mz Q_a^m \p_{Q_a^m}$, we find that the divisor equation 
still holds in the form:
\[ (p_a-zQ_a\p_{Q_a}) S^{-1}_{\tau}(mz,Q^m) = S^{-1}_{\tau}(mz,Q^m) 
(p_a\bullet_{(\tau,Q^m)}),\]
where the last subscript indicates that the matrix elements of 
$p_a\bullet$ depend on $\tau$ and $Q^m$. 
The result now follows as in Proposition 10. $\Box$

\medskip

{\tt Corollary 2.} {\em Let $\z$ be a primitive $m$th root of unity. Then
the factor $\L^{\z}=\nabla_{\z}\T^{\z}$ of the adelic cone $\hat{\L}$ is
$\D$-invariant.} 

\medskip

{\tt Proof.} Recall that the space $\T^{\z}$ is related to $\T$ by the change
$q=\z e^z$, and the action of $z$ in the operator $p_a-zQ_a\p_{Q_a}$ should be 
understood in the sense of this identification. The result follows from Lemma
since $\nabla_{\z}$ commutes with $\D$ (like $\Delta$, in Corollary 1). $\Box$   

\medskip

{\tt Theorem.} {\em Let $\D_q$ denote the algebra of finite-difference 
operators, generated by integer powers of $P_a, a=1,\dots,k$, and $Q^d$, 
with $d$ lying in the Mori cone of $X$. Define a representation of $\D_q$ 
on the symplectic loop 
space $\K$, using the operators $P_a q^{Q_a\p_{Q_a}}$    
(where $P_a$ acts by multiplication in $K^0(X)$ by the line bundle with the
Chern character $e^{-p_a}$) together with the operators of multiplication by 
$Q^d$ in the Novikov ring. Then tangent (and ruling) spaces to the overruled 
Lagrangian cone $\L\subset \K$ of true quantum K-theory on $X$ are 
$\D_q$-invariant.}
    
\medskip

{\tt Proof.} Thanks to the adelic characterization of the cone $\L$ and 
its ruling spaces, given by Theorem of Section 6 and its Corollary,
this is an immediate consequence of the following Lemma.

\medskip

{\tt Lemma.} {\em The adelic cone $\hat{L}$ is $\D_q$-invariant.} 
 
\medskip

{\tt Proof.} It is obvious that the factors $\L^{\z}$ are $\D$-invariant
for $\z$ other than roots of unity, since in this case $\L^{\z}=
\K^{\fake}_{+}$. For $\z=1$, it follows from Corollary 1 that the family of 
operators $e^{\epsilon (zQ_a\p_{Q_a}-p_a)}$ preserves $\L^{\fake}$, and so does
the operator with $\epsilon=1$, which coincides with $P_a q^{Q_a\p_{Q_a}}$.
When $\z\neq 1$ is a primitive $m$th root of unity, the family of operators
$e^{\epsilon (zQ_a\p_{Q_a}-p_a)}$ preserves $\L^{\z}$ by Corollary 2. However,
at $\epsilon=1$, the operator of the family differs from 
$P_a q^{Q_a\p_{Q_a}}$ (because $q=\z e^z$) by the factor $\z^{Q_a\p_{Q_a}}$,
which acts as $Q_a\mapsto Q_a\z$. It is essential that this extra factor
commutes with $S_{\tau^{(0)}(Q^m)}^{-1}(mz,Q^m)$ (due to $\z^m=1$). 
Since it also preserves $\H_{+}$, the result follows. $\Box$ 
 
\medskip

{\tt Example.} It is known\footnote{This result is derived from birational 
isomorphisms between some genus-0 moduli spaces of stable maps to 
$\CC P^{n-1}\times \CC P^1$ and toric compactifications of spaces of maps 
$\CC P^1\to \CC P^{n-1}$.} \cite{GiL} that for $X=\CC P^{n-1}$, 
% \J(0):=(1-q)+\sum_{a,d}\Phi_a\lan \frac{\Phi^a}{1-qL}(1-qL)\ran_{0,1}^{X,d}
\[ \J(0)=(1-q)\sum_{d=0}^{\infty} \frac{Q^d}{(1-Pq)^n\cdots (1-Pq^d)^n},\]
where $P\in K^0(\CC ^{n-1})$ represents the Hopf line bundle. 
It follows (from the string equation) that $(\J(0)/(1-q)$ lies in the
tangent space $\T_{\J(0)}\L$. Applying powers $T^r$ of the translation 
operator $T:=Pq^{Q\p_Q}$, we conclude that, for all integer $r$, the same 
tangent space contains 
\[ P^r\sum_{d=0}^{\infty} \frac{Q^dq^{rd}}{(1-Pq)^n\cdots (1-Pq^d)^n}. \]
In fact, $\J(0)$ satisfies the 2nd order finite-difference 
equation $D^n J(0)=Q\J(0)$, where $D:=1-T$. Therefore the $\D_q$-module 
generated by $J(0)/(1-q)$ is spanned over the Novikov ring by 
$T^rJ(0)/(1-q)$ with $r=0,\dots, n-1$. The projections of these elements
to $\K_{+}$ are $P^r, r=0,\dots,n-1$, which span the ring 
$K^0(\CC P^{n-1})=\ZZ [P,P^{-1}]/(1-P)^n$. 
The K-theoretic Poincare pairing on this ring is given by the residue 
formula:
\[ (\Phi(P), \Phi'(P)) = -\Res_{P=1} \frac{\Phi(P)\Phi'(P)}{(1-P)^n}
\frac{dP}{P} .\]
By computing the pairings with the above series we actually 
evaluate K-theoretic GW-invariants:
 \[ (\Phi(P), T^r\J(0)/(1-q)) = \sum_d Q^d 
\lan \frac{\Phi (P)}{1-qL}, P^r\ran_{0,2}^{X,d},\ \ r=0,\dots, n-1.\]
Thus, we started with known values of all $\lan \Phi L^k, 1\ran_{0,2}^{X,d}$ 
and computed all $\lan \Phi L^k, \Phi'\ran_{0,2}^{X,d}$ (and hence, 
by virtue of
general properties of genus-0 GW-invariants, all 
$\lan \Phi L^k, \Phi' L^l\ran_{0,2}^{X,d}$) using the $\D_q$-module 
structure alone.   
 
\section{Quantum K-theory of projective complete intersections}

{\tt Theorem.} {\em Let $X$ be a complete intersection 
in the projective space $\CC P^{n-1}$, $n>4$, given by $k(\geq 0)$ 
equations of degrees $l_1,\dots, l_k>1$, such that 
$l_1^2+\cdots l_k^2\leq n$. Then the following series represents a point in 
the overruled Lagrangian cone of true quantum K-theory of $X$:
\[ I_X:=(1-q)\sum_{d\geq 0} \frac{\prod_{j=1}^k\prod_{r=0}^{l_jd}(1-P^{l_j}q^r)}
{\prod_{r=1}^d (1-Pq^r)^n} Q^d.\]
More precisely, $I_X=\nu_*\J_X(0)$, where $\nu_*: K^0(X)\to K^0(\CC P^{n-1})$
is the K-theoretic push-forward induced by the embedding $\nu: X\to 
\CC P^{n-1}$, and $\J_X(0)$ is the value of the J-function of true quantum 
K-theory of $X$ at the input $t=0$.} 

\medskip

{\tt Remarks.} (1) To clarify this formulation,
we remind that $P$ represents the Hopf line bundle in $K^0(\CC P^{n-1})$.
%is pulled-back to $K^0(X)$, where it satisfies the unipotency relation 
%$(1-P)^{n-k}=0$. Respectively, coefficients of the $Q$-series $I_X$ lie in the 
%image of $K^0(\CC P^{n-1})$ in $K^0(X)$. 
By Lefschetz' hyperplane section 
theorem, the inclusion $X\subset \CC P^{n-1}$ induces an isomorphism
$H_2(X, \QQ) \to H_2(\CC P^{n-1},\QQ)$, whenever $2\leq n-k-2$. 
The latter holds true under our numerical restrictions on $l_j$ and $n$. 
Consequently, the degrees of holomorphic curves in $X$ are represented in 
$I_X$ by their degrees $d$ in the ambient projective space. 

(2) When $\sum l_j^2\leq n$, we also have $\sum l_j<n-2$ (strictly, unless
$k=1$, $l_1=2$, while $n=4$). Since we assumed $n>4$, we have 
for each $d>0$: 
\[ 1+\sum \frac{l_jd (l_jd+1)}{2} < n \frac{d(d+1)}{2}.\]
This means that the coefficient of $I_X$ at $Q^d$ is a reduced rational function
of $q$. Thus, the projection of $I_X$ to $\K_{+}$ is $1-q$, i.e. $I_X$ 
corresponds to the input value $t=0$. 

(3) Note that the example $n=4, k=1, l_1=2$ of the conic 
$\CC P^1\times \CC P^1 \subset \CC P^3$ is exceptional in the
sense of both previous remarks. It would be interesting to analyze the role
of the series $I_X$ in quantum K-theory of the conic.  

\medskip

{\tt Corollary.} {\em For all $s\in \ZZ$ 
 \[  \sum_d Q^d 
\lan \frac{\nu^*\Phi (P)}{1-qL}, \nu^*P^s\ran_{0,2}^{X,d} =
(\Phi(P), T^sI_X/(1-q)),\] 
where $(\cdot,\cdot)$ is the K-theoretic Poincare pairing on $K^0(\CC P^{n-1})$,
and
\[ T^s I_X/(1-q) = P^s \sum_{d\geq 0} \frac{\prod_{j=1}^k
\prod_{r=0}^{l_jd}(1-P^{l_j}q^r)}{\prod_{r=1}^d (1-Pq^r)^n}Q^dq^{sd}. \] 
}

\medskip

When $k=0$, it is known from \cite{GiL}, that the formula for $I_X$
represents the value $\J(0)$ of the K-theoretic J-function of the projective 
space. We will begin our proof of the theorem, however, with re-deriving this 
fact (and without the restriction $n>4$, of course) from the main theorem 
of this paper. After that we explain how to adjust the argument to the case
of projective complete intersections. 

To prove the theorem for $X=\CC P^{n-1}$ (let's omit the subscript $X$ in 
this case), we will show that expansions of the series $I$ near
$q=\z^{-1}$ pass the tests required by the Quantum HRR Theorem of Section 6.

Since $P$ is unipotent, $I_X$ has obviously no poles in $q$ other than roots 
of unity. 

To show that the expansion of $I$ near $q=1$ lies in the cone $\L^{\fake}$ of
the fake quantum K-theory of $X$, we begin with the following mirror-theoretic 
formula from cohomological GW-theory of projective spaces (see \cite{GiH}):
\[ \J^H(p\tau) =-ze^{-p\tau/z}\sum_{d\geq 0}
\frac{Q^d e^{d \tau}}
{\prod_{r=1}^d(p-rz)^n}.\]
Here $p$ is the hyperplane class, $\tau$ is the coordinate
on $H^2(\CC P^{n-1},\CC)$, and $\J^H$ is the J-function 
of the cohomological GW-theory of $\CC P^{n-1}$.

We now apply the ``quantum Lefschetz' theorem'' in the form of 
Coates--Corti--Iritani--Tseng \cite{CCIT} to conclude that the one-parametric 
family  
\[ I^{\fake}(\tau):=-z P^{-\tau/z}\sum_{d\geq 0} \frac{Q^d e^{d\tau}}
{\prod_{r=1}^d(1-Pe^{rz})^n},\ \ z=\log q, \ P=\ch^{-1}(e^{-p}),\]
lies in $\L^{\fake}$, the overruled Lagrangian cone of fake quantum K-theory
of $\CC P^{n-1}$. 

Let us recall from Section 5 that $\L^{\fake}=\ch^{-1}\Delta \L^H$, where
\[ \log \Delta \sim \sum_{r=1}^{\infty}\sum_xs(x-rz), \ \ 
s(u):=\log \frac{u}{1-e^{-u}},\]
and $x$ runs Chern roots of $T_{\CC P^{n-1}}$. We claim that in fact 
\[ \log \Delta \sim \sum_{r=1}^{\infty} \left( n\ s(p-rz) - s(-rz)\right).\]
Indeed, $T_{\CC P^{n-1}} =nP^{-1}-1$, 
and the construction 
of the operator $\log \Delta$ from (Chern roots of) a bundle is additive.
Note that the last summand does not affect the way $\Delta$ acts on $\L^H$,
since, being overruled, the cone $\L^H$ is invariant under multiplication by 
functions of $z$.   

Replace (following \cite{CCIT}) the degree-2 cohomology class $p$ with 
the operator of differentiation $-z\p_p$ in the direction of this class, and
consider the Euler--Maclaurin asymptotics: 
\[ \log \hat{\Delta} \sim\sum_{r=1}^{\infty} n\ s(-z\p_p-rz).\] 
On the one hand, applying $\hat{\Delta}^{-1}$ to the J-function $\J^H$ (i.e.
making $\p_p$ act as the derivative in $\tau$, one obtains a family still 
lying in the same cone $\L^H$ as $\J^H$.\footnote{This is the result of 
a non-trivial lemma from \cite{CGi} based on properties of the D-module 
generated by the J-function and used there in the proof of a Quantum
Lefschetz' Theorem, generalized in \cite{CCIT}.} On the other hand, 
\[ s(-z\p_p-rz) e^{d\tau-p\tau/z}=s(p\tau-(r+d)z) e^{d\tau-p\tau/z},\]
and hence 
\[ \hat{\Delta}^{-1} \left(e^{d\tau-p\tau/z}\right) = \Delta^{-1} \times 
\left(e^{d\tau-p\tau/z}\right) \prod_{r=1}^d e^{n\ s(p-rz)}.\]
It follows that modifying each $Q^d$-term in $\J^H$ by the factor
\[ \prod_{r=1}^de^{n\ s(p-rz)} = \prod_{r=1}^d \frac{(p-rz)^n}{(1-e^{-p+rz})^n},\]
one obtains a family lying in the cone $\Delta \L^H = \ch(\L^{\fake})$.
The modified family is $I^{\fake}(\tau)$ indeed.

Note that $(1-e^z) I^{\fake}(0)/(-z)$ coincides with $I$ after the change 
$e^{-p}=:P$ and $e^z=:q$, that is, {\em the expansion of $I$ near $q=1$ 
coincides with $ (1-e^z) I^{\fake}(0)/(-z) \in \L^{\fake}$} as required.  
        
\medskip

To analyze the expansion of $I$ near $q=\z^{-1}$ where $\z$ is an $m$-th root of
$1$, we note first, that by the previous the family $I^{\fake}(\tau)/(-z)$
represents tangent vectors $\L^{\fake}$ at the family of application points 
$(1-q)I^{\fake}(\tau)/(-z)$. Following test (iii) in Theorem of Section 6, 
we first change $q\mapsto q^m$, $Q\mapsto Q^m$, and obtain the following 
family of tangent vectors to the cone thus transformed: 
\[ e^{-p\tau/\log q^m}\sum_{d\geq 0} \frac{Q^{md}e^{d\tau}}{(1-Pq^m)^n(1-Pq^{2m})^n\cdots 
(1-Pq^{md})^n}.\]
Now we employ the construction of \cite{CCIT} once again, this time using the 
operator $\nabla_{\z}$. Namely, we put:
\[ \hat{\nabla}_{\z} \sim_{q\z\to 1} \prod_{r=1}^{\infty} \frac{(1-q^{mr+m\p_p})^n}
{(1-q^{r+m\p_{p}})^n},\] 
and apply $\hat{\nabla}_{\z}^{-1}$ to the above family of tangent vectors.
Since 
\[ q^{m\p_p} \left(e^{-p\tau/\log q^m}e^{d\tau}\right) = e^{-p\tau/\log q^m}e^{d\tau} 
P q^{md},\]
and hence
\[ \hat{\nabla}_{\z}^{-1} \left(e^{-p\tau/\log q^m}e^{d\tau}\right) = 
\nabla_{\z}^{-1}\times \left(e^{-p\tau/\log q^m}e^{d\tau}\right)
 \frac{\prod_{r=1}^d(1-Pq^{mr})^n}{\prod_{r=1}^{md}(1-Pq^r)^n}.\]
The result is $\nabla_{\z}^{-1}\times \tilde{I}_{\z}(\tau)$, where 
\[ \tilde{I}_{\z}:= e^{-p\tau/\log q^m} \sum_{d\geq 0} \frac{Q^{md} e^{d\tau}}
{(1-Pq)^n(1-Pq^2)^n\cdots (1-Pq^{md})^n}.\] 
The expression should be interpreted as a Laurent series expansion
near $q=\z^{-1}$, and as such should be compared with $I_{\z}$,
where 
\[ I = \sum_{d\geq 0} \frac{Q^d}{(1-Pq)^n(1-Pq^2)^n\cdots (1-Pq^d)^n}.\]
At $\tau=0$, the terms of $\tilde{I}_{\z}$ reproduce correctly
the terms of $I_{\z}$ with $d$ divisible by $m$, but all other terms are 
missing. 

Nevertheless, we derive from this that $I_{\z}$ lies in the space 
$\nabla_{\z} \T$ 
where $\T$ is obtained from $T_{I^{\fake}(0)}\L^{\fake}$ by the changes 
$q\mapsto q^m, Q\mapsto Q^m$. Namely, let $\T_{\tau}$ denote the spaces thus 
obtained from $T_{I^{\fake}(\tau)}\L^{\fake}$. We have: 
$\tilde{I}_{\z}(\tau) \in \nabla_{\z}\T_{\tau}$. Consider operators
\[ D(\tau):=\sum_{d=0}^{m-1}\frac{Q^de^{d\tau}}{\prod_{r=1}^d(1-q^{-m\p_p}q^r)^n}.\]
It should be understood as an expansion near $q=\z^{-1}$, and it is important 
that within the given range $0<r\leq d<m$ of the indices $d$ and $r$, 
the denominators have no zeroes at $q=\z^{-1}$, and thus $D(\tau)$ is a 
power series in $z\p_p$ ($z=\log q$). We conclude that  
$D(\tau) \tilde{I}_{\z}(\tau) \in \nabla_{\z}\T_{\tau}$ (since $z\p_p$
preserves tangent spaces to overruled cones). At $\tau=0$, we have:  
\[ I_{\z} = D(\tau) \tilde{I}_{\z}(\tau)\ |_{\tau=0},\]
which thus proves that $I_{\z} \in \nabla_{\z} \T$.
   
\medskip

What we have established about the series $I$ means that the decomposition of 
it into elementary fractions obeys the recursion relations of Section 7, 
with the leg contribution obtained by Adams' operation $\Psi$ from the arm
contribution, corresponding to the input point $t(q)=[I]_{+}-(1-q)$. Since 
the projection $[\dots]_{+}$ of $I$ to $\K_{+}$ is $1-q$, we find that $t=0$
as required, and hence $I=\J(0)$. $\square$ 

\medskip

{\tt Remark.} With the exception of the last property $[I]_{+}=1-q$, this
seemingly sophisticated argument is in fact general enough to work for 
$q$-hypergeometric series $I_X$ that can be associated to any symplectic toric 
manifold $X$ as follows. Let $X$ be obtained by symplectic reduction
$X=\CC^n//T^k$ by the action of the subtorus $T^k\subset T^n$ of the maximal 
torus, the embedding being determined (in some basis of $\pi_1(T^k)$) by the 
integer matrix $(m_{ij})$, $i=1,\dots,k$, $j=1,\dots, n$ (see \cite{GiZ, GiT} 
for more details). Let $Q^d=Q_1^{d_1}\cdots Q_k^{d_k}$ represents a point in the 
Mori cone of $X$ in coordinates $(d_1,\dots,d_k)$ on $H_2(X)$ corresponding to
the chosen basis of $\pi_1(T^k)$, and  $P_i^{-1}$, $i=1,\dots,k$, denote the 
line bundles over $X$ whose 1st Chern classes form the dual basis of $H^2(X)$. 
In this notation: 
\[ I_X =\sum_{d} Q^d\prod_{j=1}^n\frac{\prod_{r=-\infty}^{0}
(1-q^r\prod_{i=1}^kP_i^{m_{ij}})}{\prod_{r=-\infty}^{\sum_id_im_{ij}}
(1-q^r\prod_{i=1}^kP_i^{m_{ij}})}.\]    
The property $[I_X]_{+}=1-q$, however, does not hold unless $X$ is a product of
complex projective spaces. It would be interesting to find out if nevertheless
$I_X\in \L_X$. 

\medskip

The above computation will also work for the series $I_X$ corresponding to
projective complete intersection described in the theorem. However, there is 
a catch here, related to the fact that cohomology and K-theory of $X$ may not
be entirely describable in terms of the ambient projective space, and thus the
information gained about $I_X$ won't yet allow to make a legitimate application
of our Quantum HRR Theorem. More specifically, our computation would only 
be concerned with the properties of $\nu_*(I_X)$ expressed in terms of 
$\nu_*(I^{\fake})$, and the latter may not even lie on $\L^{\fake}$.

In order to bypass the difficulty, we introduce a model of quantum K-theory 
of a {\em supermanifold} $\Pi E$, interpolating between those of $X$ and 
$\CC P^{n-1}$. Let $E$ be the total space of the sum of the line bundles 
over $\CC P^{n-1}$ of degrees $l_1,\dots,l_k$, while $\Pi$ indicates the 
fiberwise parity change. By definition, genus-0 moduli spaces of stable maps
to $\Pi E$ are the same as to $\CC P^{n-1}$, but the virtual structure sheaf
is changed, by tensoring the structure sheaf ${\mathcal O}^{vir}_{0,r,d}$ 
with {\em the $S^1$-equivariant K-theoretic Euler class
of the bundle $E_{0,r,d}$ (i.e. the Koszul complex of the dual, $E^*_{0,r,d}$).} 
Here $E_{0,r,d}$ stands for the bundle $\pi_*\ev^*E$ whose fiber over a stable map
$f: \Sigma \to \CC P^{n-1}$ is $H^0(\Sigma, f^*E)$. The circle $S^1$ is made to
act by multiplication by unitary scalars fiberwise on $E$, and hence on 
$E_{0,r,d}$. Respectively, correlators of quantum K-theory of $\Pi E$ take 
values in the representation ring $\CC [S^1] = \CC [\Lambda, \Lambda^{-1}]$.
Their algebraic-geometrical meaning (instead of holomorphic Euler 
characteristics of a sheaf) is the trace of $S^1$ on the sheaf cohomology. 
The ring $K^0(\Pi E)$ coincides with $K^0(\CC P^{n-1}) \otimes 
\CC [S^1]$, and is equipped with the K-theoretic Poincare pairing
\[ (\Phi ,\Phi')_{\Pi E} = -\Res_{P=1} \Phi(P)\Phi'(P) \frac{\prod_{j=1}^k 
(1-P^{l_j}\Lambda)}{(1-P)^n}\frac{dP}{P}.\]
This pairing becomes non-degenerate if division by $1-\Lambda$ is allowed. 
After this localization, the resulting quantum K-theory of the supermanifold 
$\Pi E$ satisfies all the axioms of genus-0 quantum K-theory. 

Furthermore, the Quantum HRR Theorem of Section 6 and its proof given in 
Sections 7 and 8 work verbatim for true quantum K-theory of 
$\Pi E$.\footnote{Note that we are not using any geometric fixed point
localization with respect to $S^1$, so that all moduli spaces, Kawasaki strata,
etc. remain the same, and only the meaning and values of the correlators are
modified appropriately.} 

Thus, applying the same technology as we did in the case of $X=\CC P^{n-1}$, we 
establish that {\em under the numerical assumptions of Theorem, 
$\J_{\Pi E}(0)=I_{\Pi E}$, where} 
\[ I_{\Pi E}:=\sum_{d\geq 0} Q^d \frac{\prod_{j=1}^k \prod_{r=1}^{l_jd}
(1-P^{l_j}\Lambda q^r)}{\prod_{r=1}^d(1-Pq^r)^n}.\]
Here are some formulas that elucidate this claim:
\[ \J^H_{\Pi E}(p\tau) =-ze^{-p\tau/z}\sum_{d\geq 0}
Q^d e^{d \tau} \frac{\prod_{j=1}^k\prod_{r=1}^{l_jd}(\lambda+l_jp-rz)}
{\prod_{r=1}^d(p-rz)^n},\]
where $\lambda$ is the 1st Chern class of the universal $S^1$-bundle 
$\Lambda^{-1}$ (i.e. $\ch (\Lambda) = e^{-\lambda}$);
\[ \hat{\Delta} \sim \prod_{r=1}^{\infty}\frac{(-z\p_p-rz)^n}{(1-e^{z\p_p+rz})^n}
\prod_{j=1}^k\frac{(1-e^{-\lambda+l_jz\p_p+rz})}{(\lambda-l_jz\p_p-rz)};\]
\[ I^{\fake}_{\Pi E}(\tau):=-(\log q) P^{-\tau/\log q}\sum_{d\geq 0} Q^d e^{d\tau}
\frac{\prod_{j=1}^k\prod_{r=1}^{l_jd}(1-P^{l_j}\Lambda q^r)}{\prod_{r=1}^d(1-Pq^r)^n}
;\]  
\[ \hat{\nabla}_{\z} \sim_{q\z\to 1} \prod_{r=1}^{\infty} \frac{(1-q^{mr+m\p_p})^n}
{(1-q^{r+m\p_{p}})^n} \prod_{j=1}^k\frac{(1-\Lambda q^{r+ml_j\p_p})}
{(1-\Lambda q^{mr+ml_j\p_p})};\] 
\[ \tilde{I}_{\z}:= e^{-p\tau/\log q^m} \sum_{d\geq 0} Q^{md} e^{d\tau}
\frac{\prod_{j=1}^k\prod_{r=1}^{ml_jd}(1-\Lambda P^{l_j}q^r)}
{\prod_{r=1}^{md}(1-Pq^r)^n}; \] 
\[ D(\tau):=\sum_{d=0}^{m-1}Q^de^{d\tau}\frac{\prod_{j=1}^k\prod_{r=1}^{l_jd}
(1-\Lambda q^{-ml_j\p_p}q^r)}{\prod_{r=1}^d(1-q^{-m\p_p}q^r)^n}.\]

\medskip

Once the equality $\J_{\Pi E}(0)=I_{\Pi E}$ is proved, to establish the equality 
$\nu_* \J_X(0)=I_X$, it remains to notice that for all $s\in \ZZ$
\[ (\nu^* P^s, \J_X(0))_X = (P^s, \J_{\Pi E}(0))_{\Pi E}\left|_{\Lambda=1}. \right. \]
Indeed, when $X$ is given in $\CC P^{n-1}$ by a section of $E$, the moduli 
space $X_{0,r,d}$ is given in $(\CC P^{n-1})_{0,r,d}$ by the corresponding section
of the bundle $E_{0,r,d}$, and (according to \cite{To, To1}) the virtual 
structure sheaf of $X_{0,r,d}$ is described in $K^0((\CC P^{n-1})_{0,r,d})$ 
by tensoring the virtual structure sheaf of $(\CC P^{n-1})_{0,r,d}$ with 
the K-theoretic Euler class of $E_{0,r,d}$, albeit, the non-equivariant one, and 
hence the specialization to $\Lambda=1$. 

\medskip

{\tt Remark.} In a separate paper, we will describe a different machinery,
based on the Adams--Riemann--Roch--Gr\"othendieck formula, 
that allows one to compare genus-0 
K-theoretic GW-invariants of complete intersections (and more general 
``twisted'' quantum K-theories) with those of the ambient space, 
remaining entirely on K-theoretic grounds.

\enddocument